\documentclass{article}
\usepackage{amsmath,amssymb,amsthm,graphicx,epsfig,subfigure,float,url}
\usepackage[colorlinks=true]{hyperref}
\usepackage{pdfsync}
\usepackage{color}

\usepackage{verbatim}
\usepackage[applemac]{inputenc}
 
\usepackage[usenames,dvipsnames]{pstricks}
\usepackage{pst-grad} 
\usepackage{pst-plot} 

\def\R{\textrm{I\kern-0.21emR}}
\def\N{\textrm{I\kern-0.21emN}}

\newcommand{\C} {\mathbb{C}}

\renewcommand{\geq}{\geqslant}
\renewcommand{\leq}{\leqslant}

\newtheorem{theorem}{Theorem}

\newtheorem{proposition}{Proposition}

\newtheorem{lemma}{Lemma}
\theoremstyle{definition}
\theoremstyle{definition}\newtheorem{remark}{Remark}

\newcommand{\Hun}{\mathbf{(H_1)}}
\newcommand{\Hdeux}{\mathbf{(H_2)}}
\newcommand{\Htrois}{\mathbf{(H_3)}}
\newcommand{\Real}{\mathrm{Re}}

\title{Actuator design for parabolic distributed parameter systems with the moment method}

\author{Yannick Privat\footnote{CNRS, Sorbonne Universit\'es, UPMC Univ Paris 06, UMR 7598, Laboratoire Jacques-Louis Lions, F-75005, Paris, France ({\tt yannick.privat@upmc.fr}).}
	\and Emmanuel Tr\'elat\footnote{Sorbonne Universit\'es, UPMC Univ Paris 06, CNRS UMR 7598, Laboratoire Jacques-Louis Lions, Institut Universitaire de France, F-75005, Paris, France (\texttt{emmanuel.trelat@upmc.fr}).} 
        \and Enrique Zuazua\footnote{DeustoTech - Fundaci\'on Deusto, Avda Universidades, 24, 48007, Bilbao - Basque Country - Spain.}~\footnote{Departamento de Matem\'aticas, Universidad Aut\'onoma de Madrid, 28049 Madrid - Spain}~\footnote{Facultad Ingenieria, Universidad de Deusto, Avda. Universidades, 24, 48007, - Basque Country - Spain.}~\footnote{Sorbonne Universit\'es, UPMC Univ Paris 06, CNRS UMR 7598, Laboratoire Jacques-Louis Lions, F-75005, Paris, France (\texttt{enrique.zuazua@uam.es}).}
}

\date{}

\begin{document}

\maketitle

\begin{abstract}
In this paper, we model and solve the problem of designing in an optimal way actuators for parabolic partial differential equations settled on a bounded open connected subset $\Omega$ of $\R^n$. 
We optimize not only the location but also the shape of actuators, by finding what is the optimal distribution of actuators in $\Omega$, over all possible such distributions of a given measure.
Using the moment method, we formulate a spectral optimal design problem, which consists of maximizing a criterion corresponding to an average over random initial data of the largest $L^2$-energy of controllers. Since we choose the moment method to control the PDE, our study mainly covers one-dimensional parabolic operators, but we also provide several examples in higher dimensions.

We consider two types of controllers: either internal controls, modeled by characteristic functions, or lumped controls, that are tensorized functions in time and space.
Under appropriate spectral assumptions, we prove existence and uniqueness of an optimal actuator distribution, and we provide a simple computation procedure. Numerical simulations illustrate our results.
\end{abstract}

\noindent\textbf{Keywords:} heat equation, parabolic systems, shape optimization, null controllability, moment method, lumped control.

\medskip

\noindent\textbf{AMS classification:} 93B07, 35L05, 49K20, 42B37.
\section{Introduction and modeling of the problem}\label{secintro}

In this article, we model and solve the problem of finding the optimal shape and location of internal controllers for parabolic equations with (mainly) Dirichlet boundary conditions and (mainly) in the one-dimensional case $\Omega=(0,\pi)$. Such questions are frequently encountered in engineering applications. We provide a possible mathematical model for investigating such issues.

For mathematical reasons that will be clarified in the sequel, we will focus in the whole article on  controls obtained by using the so-called {\it moment method}. As it will be underlined, it requires in general some spectral gap assumptions on the operators involved that essentially reduce the applications of our results to one-dimensional partial differential equations, but our results also cover several particular situations in larger dimension.

To avoid technicalities and highlight the main ideas, we first present the results in the simplified framework of the controlled one-dimensional heat equation with Dirichlet boundary conditions, without introducing (at this step) the more general parabolic framework in which our results are actually valid. 

Generalizations to a more general framework will be described in Section \ref{sec:Control}. 
Unlike the simplified case of the one-dimensional heat equation, it requires a discussion on the M\"untz-Sz\'asz theorem as well as specific spectral considerations.

Notice that the general control framework in which this problem could be addressed is much more intricate and will be evoked as a possible perspective at the end of this article. 

\subsection{Reminders on the controllability of the 1D heat equation}\label{sec:control2001}
Consider the internally controlled one-dimensional heat equation
\begin{equation}\label{heatEqcontrolled_intro}
\partial_t y(t,x)-\partial_{xx} y(t,x)=\chi_\omega (x)u(t,x), \quad (t,x)\in (0,T)\times (0,\pi),
\end{equation}
with Dirichlet boundary conditions
\begin{equation}\label{dir_cond}
y(t,0)=y(t,\pi)=0, \qquad t\in (0,T),
\end{equation}
where $u\in L^2((0,T)\times(0,\pi))$ is a control function, and $\omega$ is a measurable subset of $(0,\pi)$ standing for the support of the controller. Here, $\chi_\omega$ is the characteristic function of $\omega$, defined by $\chi_\omega(x)=1$ if $x\in\omega$ and $\chi_\omega(x)=0$ otherwise. For a given subset $\omega$, the equation \eqref{heatEqcontrolled_intro} is said to be exactly null controllable in time $T$ whenever every initial datum $y(0,\cdot)\in L^2(0,\pi)$ can be steered to $0$ in time $T$ by means of an appropriate control function $u\in L^2((0,T)\times(0,\pi))$. It is well known that, for a given subset $\omega$, the system \eqref{heatEqcontrolled_intro} is exactly null controllable if and only if there exist a positive constant $C$ (only depending on $T$ and $\omega$) such that
\begin{equation}\label{ineqobs}
C \int_0^\pi z(T,x)^2\, dx
\leq \int_0^T\int_\omega z(t,x)^2 \,dx \, dt,
\end{equation}
(observability inequality)
for every solution of 
\begin{equation}\label{heatEq}
\begin{split}
& \partial_t z(t,x)-\partial_{xx} z(t,x)=0, \quad (t,x)\in (0,T)\times(0,\pi),\\
& z(t,0)=z(t,\pi)=0, \qquad\qquad t\in (0,T),
\end{split}
\end{equation}
with $z(0,\cdot)\in L^2(0,\pi)$. 

\paragraph{Exact null controllability by the moment method.}
The observability inequality \eqref{ineqobs} has been shown to hold true for any subset $\omega$ of $[0,\pi]$ of positive Lebesgue measure in \cite{Russell} by the moment method, that we will use as well in the present paper and that we recall hereafterin.

The eigenfunctions of the Dirichlet-Laplacian, given by $\phi_j(x)=\sqrt{\frac{2}{\pi}}\sin(jx)$ for every $j\in\N^*$, associated with the eigenvalues $\lambda_j=j^2$, make up an orthonormal basis of $L^2(0,\pi)$. From the M\"untz-Sz\'asz theorem, there exists a sequence $(\theta_j^T)_{j\in\N^*}$ of $L^2(0,T)$, biorthogonal to the sequence of functions $t\mapsto e^{-j^2 t }$.
The following lemma provides an exact null controllability result for \eqref{heatEqcontrolled_intro}-\eqref{dir_cond}.

\begin{lemma}\label{lemm:contHeat}\cite{Russell}
Let $T>0$ and let $\omega$ be a measurable subset of $(0,\pi)$ of positive measure. Then every initial datum
$$
y(0,x)=y^0(x)=\sum_{j=1}^{+\infty} a_j\sin(jx),
$$
in $L^2(0,\pi)$, can be steered to zero in time $T$ with the control $u\in L^2((0,T)\times (0,\pi))$ defined by
$$
u(t,x)=-\sum_{j=1}^{+\infty} \frac{a_j e^{-j^2T}}{\int_\omega\sin^2(jy) \, dy}\theta_j^T(T-t)\sin(jx) .
$$
\end{lemma}

A proof of this lemma is given in a more general setting in Section \ref{secproof31}.

We set $\Gamma_{\omega}(y^0)=\chi_\omega u$. The operator $\Gamma_{\omega}:L^2(0,\pi)\rightarrow L^2((0,T)\times(0,\pi))$ is linear and continuous, and is called the \textit{moment control operator}.
The norm of this operator, given by
$\Vert \Gamma_{\omega}\Vert = \sup \{ \Vert \Gamma_{\omega} (y^0)\Vert_{L^2((0,T)\times(0,\pi))} \mid \Vert y^0\Vert_{L^2(0,\pi)}=1 \},$
provides an account for the worst possible initial datum to be controlled to zero, in terms of the effort ($L^2$ energy) required to steer this initial datum to zero.
Minimizing $\Vert \Gamma_{\omega}\Vert$ over a class of admissible domains (that we will denote by $\mathcal{U}_L$ in the sequel) is then an interesting problem, that will be discussed in the next section.

\subsection{State of the art}
When realizing exact null controllability in practice, an important question is to know where to place and how to shape optimally the actuators (modeled here by the subset $\omega$), in order to minimize the efforts done to steer any possible initial data to zero. In this paper, we want to optimize not only the location, but also the shape of actuators, without any specific restriction on the regularity of $\omega$.

The literature on optimal sensor or actuator location problems is abundant in engineering applications (see, e.g., \cite{armaoua,harris,Kumar,Morris,Sigmund,Ucinski,vandeWal,vandewouwer} and references therein), where the aim is often to optimize the number, the place and the type of sensors or actuators in order to improve the estimation of the state of the system. Fields of applications are very numerous and concern for example active structural acoustics, piezoelectric actuators, vibration control in mechanical structures, damage detection and chemical reactions, just to name a few of them. In most of these applications the method consists of approximating appropriately the problem by selecting a finite number of possible optimal candidates and of recasting the problem as a finite-dimensional combinatorial optimization problem. In many of these contributions the sensors or actuators have a prescribed shape (for instance, balls with a prescribed radius) and then the problem consists of placing optimally a finite number of points (the centers of the balls) and thus is finite-dimensional, since the class of optimal designs is replaced with a compact finite-dimensional set. We stress that, in the present paper, the shape of the control domain is an unknown of the optimization procedure.

From the mathematical point of view, the issue of studying a relaxed version of optimal design problems for the shape and position of sensors or actuators has been investigated in a series of articles. In \cite{munchHeat}, the authors study a homogenized version of the optimal location of controllers for the heat equation problem (for fixed initial data), noticing that such problems are often ill-posed. In \cite{allaireMunch}, the authors consider a similar problem and study the asymptotic behavior as the final time $T$ goes to infinity of the solutions of the relaxed problem; they prove that optimal designs converge to an optimal relaxed design of the corresponding two-phase optimization problem for the stationary heat equation. We also mention \cite{munchPedr} where, for fixed initial data, numerical investigations are used to provide evidence that the optimal location of null-controllers of the heat equation problem is an ill-posed problem.

Concerning the problem of optimal shape and location of sensors for fixed initial data (instead of controllers in \cite{munchPedr}) we proved in \cite{PTZobspb1} that it is always well posed for heat, wave or Schr\"odinger equations (in the sense that no relaxation phenomenon occurs); we showed that the complexity of the optimal set depends on the regularity of the initial data, and in particular we proved that, even for smooth initial data, the optimal set may be of fractal type (and there is no relaxation). In \cite{PTZparabND}, we modeled and solved the problem of optimal shape and location of the observation domain having a prescribed measure. This problem was motivated by the question of shaping and placing sensors in some domain in such a way to optimize the quality of the observation.

Here, we rather investigate the dual question of the best shape and location of actuators.
In \cite{munchZuazua}, the authors investigate numerical approximations of exact or trajectory controls for the heat equation, by developing a numerical version of the so-called transmutation method.

\subsection{Modeling of the optimal design problems: a randomization procedure}
In the present paper, our objective is to search the internal control domain over \textit{all} possible subsets of $(0,\pi)$, without assuming any a priori regularity. We optimize not only the placement but also the shape of the actuators.

\medskip

Note that, for any problem consisting of optimizing the quality of the control, certainly the best strategy consists of controlling the solutions over the whole domain $(0,\pi)$. This is however obviously not reasonable and in practice the domain covered by actuators is limited, due for instance to cost considerations. From the mathematical point of view, we model this basic limitation by considering as the set of unknowns, the set of all possible measurable subsets $\omega$ of $(0,\pi)$ that are of Lebesgue measure $\vert\omega\vert = L\pi$, where $L\in(0,1)$ is some fixed real number. Any such subset $\omega$ represents the actuators put in $(0,\pi)$. Finally, for mathematical reasons, it is more convenient to assimilate a measurable subdomain $\omega$ of $(0,\pi)$ to its characteristic function $\chi_{\omega}$, vanishing outside $\omega$ and equal to 1 else. Hence, let us introduce the class of admissible control domains
\begin{equation}\label{defUL}
\boxed{
\mathcal{U}_L = \{ \chi_\omega\in L^\infty(0,\pi,\{0,1\})\ \vert\ \omega\subset \Omega\ \textrm{measurable} ,\  \vert\omega\vert=L \pi\} .
}
\end{equation}

In view of modeling the optimal design of actuators, a first approach consists of minimizing the functional 
$\chi_{\omega}\mapsto \Vert \Gamma_{\omega}\Vert$
over the set $\mathcal{U}_{L}$. However, even for simple choices of control domains $\omega$, the quantity $\Vert \Gamma_{\omega}\Vert$ is not explicitly computable and therefore the cost functional is hard to handle. Besides, note that the moment control operator norm $\Vert \Gamma_{\omega}\Vert$ is \textit{deterministic} and thus provides an account for the \textit{worst possible case}; in this sense, it is a \textit{pessimistic} constant. One can argue that, in practice, when running a large number of experiments, it is expected that the worst possible case does not occur so often. For these reasons, we are next going to consider an average criterion which, in some sense, does not take into account rare events. Nevertheless, we stress that the issue of minimizing $\Vert \Gamma_{\omega}\Vert$ with respect to the domain $\omega$ has not only a mathematical interest, but appears also naturally in some practical situations, where it is imperative that the worst possible case be avoided, even if it is a rare event. We refer to the conclusion section \ref{sec:conclpb} for some comments about such a problem.
The same kind of difficulty arises when modeling optimal design problems for sensors, as discussed in \cite{PTZobsND,PTZparabND}.

In this paper, we propose another approach based on the controllability result stated in Lemma \ref{lemm:contHeat}, and on a randomization argument reflecting what happens when a large number of experiments is expected to be done.
We are going to use a probabilistic argument, by considering random initial data. We follow the approach developed in \cite{PTZobsND,PTZparabND}. 
Let us fix an arbitrary initial datum $y(0,\cdot)=y^0(\cdot)\in L^2(0,\pi)$ of the controlled system \eqref{heatEqcontrolled_intro}, with Fourier coefficients defined by 
\begin{equation}\label{defajbj}
a_j = \int_0^\pi y^0(x) \sin (jx)\, dx,
\end{equation}
These coefficients are now randomized according to $a_j^\nu=\beta^\nu_{j}a_j$ for every $j\in\N^*$, where $(\beta_{j}^\nu)_{j\in\N^*}$ is a sequence of independent real random variables on a probability space $(\mathcal{X},\mathcal{F},\mathbb{P})$, having mean equal to $0$, variance equal to $1$, and with a super exponential decay\footnote{Recall that the sequence $(\beta_{j}^\nu)_{j\in\N^*}$ is said to have a {\it super-exponential decay} whenever
$$
\exists (C, \delta)\in (\R_+^*)^2\ \mid \ \forall \alpha \in \R, \ \mathbb{E}(e^{\alpha |\beta_{j}^\nu|})\leq Ce^{\delta \alpha^2}.
$$} (for instance, independent Bernoulli random variables, see \cite{Burq,BurqTzvetkov1} for  details and properties of randomization).
For every event $\nu\in\mathcal{X}$, the control steering the initial datum 
$$
y^0_\nu (x)= \sum_{j=1}^{+\infty}\beta_j^\nu a_j \sin(jx)
$$
to zero by the moment method is, according to Lemma \ref{lemm:contHeat},
$$
u^\nu (t,x)=-\sum_{j=1}^{+\infty} \frac{\beta_j^\nu a_j e^{-j^2T}}{\int_\omega\sin^2(jy) \, dy}\theta_j^T(T-t)\sin(jx).
$$
Using the previous notations, one has $\Gamma_{\omega}(y^0_\nu )=\chi_\omega u^\nu$. We propose, then, to model the problem of best actuator shape and location as the problem of minimizing the averaged functional
$$
\mathcal{K}(\chi_\omega) = \sup_{\Vert y^0\Vert_{L^2(0,\pi)}=1} \mathbb{E}\left(\Vert \Gamma_{\omega}(y^0_{\nu})\Vert^2_{L^2((0,T)\times (0,\pi))}\right)
$$
over $\mathcal{U}_L$, where $\mathbb{E}$ is the expectation over the probability space $(\mathcal{X},\mathcal{F},\mathbb{P})$.
This is the randomized counterpart to the deterministic quantity $\Vert \Gamma_{\omega}\Vert$.
One of the advantages is that $\mathcal{K}(\chi_\omega)$ can be explicitly computed, as follows.

\begin{lemma}\label{lem2}
One has
$$
\mathcal{K}(\chi_\omega) =\left( \inf_{j\in\N^*} \frac{e^{2j^2T}}{\int_{0}^T\theta_j^T(t)^2\, dt} \int_{\omega}\sin^2(jx)\, dx\right)^{-1} ,
$$
for every measurable subset $\omega\subset(0,\pi)$.
\end{lemma}

Lemma \ref{lem2} is proved in Section \ref{sec_proof_lem2}.
Therefore, the problem of best shape and location of actuators is finally written as
\begin{equation}\label{def:odp1D}
\boxed{
\sup_{\chi_{\omega}\in\mathcal{U}_{L}}\inf_{j\in\N^*} \frac{e^{2j^2T}}{\int_{0}^T\theta_j^T(t)^2\, dt} \int_{\omega}\sin^2(jx)\, dx.
}
\end{equation}

The article is organized as follows: Section \ref{sec31} is devoted to comments on the control of parabolic equations by the moment method, the use of biorthogonal sequences and modeling of the optimal design problem issues. In Section \ref{sec32}, we solve the problem and provide a numerical illustration as well as a series of examples, mainly in the 1D case due to the restrictions imposed by the choice of the control method.
Finally, in Section \ref{sec:OLC}, we investigate a variant of the previously studied optimal design problem, where the control acts on the system by means of tuning the time-intensity.

\section{Solving the problem \eqref{def:odp1D}}\label{sec:mainResIntro}
\subsection{Main results, comments and illustration}\label{sec:mainrescomillu}
We first provide an existence result.

\begin{theorem}\label{thm:odp1D}
The shape optimization problem \eqref{def:odp1D} has a unique\footnote{Here and in the sequel, it is understood that the optimal set is unique within the class of all measurable subsets of $(0,\pi)$ quotiented by the set of all measurable subsets of $\Omega$ of zero measure.} solution $\chi_{\omega^*}$. 
\end{theorem}

This theorem is proved in Section \ref{sec_thm:odp1D}.

In addition to this result, what is remarkable is that we have a simple and numerically efficient procedure to compute the optimal control domain $\omega^*$.

\paragraph{Algorithmic computation procedure.} 
The optimal set $\omega^*$ of Theorem \ref{thm:odp1D} can actually be built from a finite-dimensional spectral approximation, by keeping only a finite number of modes. Let us provide the details of the procedure.
For every integer $N\in\N^*$, we define the functional $\mathcal{J}_N$ by
\begin{equation*}
\mathcal{J}_N(\chi_\omega) = \inf_{1\leq j\leq N} \frac{e^{2j^2T}}{\int_{0}^T\theta_j^T(t)^2\, dt} \int_{\omega} \sin^2(jx)\, dx,
\end{equation*}
for every measurable subset $\omega$ of $(0,\pi)$.
The functional $\mathcal{J}_N$ is a spectral truncation to the $N$ first terms.
We consider the shape optimization problem
\begin{equation}\label{pb_max_JN}
\sup_{\chi_\omega\in {\mathcal{U}}_L} \mathcal{J}_N(\chi_\omega),
\end{equation}
called truncated problem, which is a spectral approximation of the problem \eqref{def:odp1D}.

We have then the following results, proved in Sections \ref{sec_A2} and \ref{sec_thm:odp1D}.

\begin{proposition}\label{truncTheo}
For every $N\in\N^*$, the truncated problem \eqref{pb_max_JN} has a unique solution $\chi_{\omega^N}\in\mathcal{U}_L$.
Moreover, $\omega^N$ has a finite number of connected components, and there exists $\eta^N>0$ such that $\omega^N\subset [\eta^N,\pi-\eta^N]$.
\end{proposition}

Proposition \ref{prop:trunc} further (see Section \ref{sec32}) will provide an extension of this result to higher dimensions. We will however provide two different proofs. Indeed, in the one-dimensional case investigated here, we will show in the proof that the problem \eqref{pb_max_JN} can be expressed in an equivalent way as a classical optimal control problem. This point of view (already used in \cite{PTZ_HUM1D}) is interesting not only for the proof but also in order to derive efficient numerical methods for the numerical computation of the optimal domains.

Let us now give the main result that is at the base of the algorithmic procedure.

\begin{theorem}\label{thm:stat1D}
There exists $N_0\in \N^*$ such that $\omega^*=\omega^N$ for every $N\geq N_0$.

Furthermore, we have $N_0\leq \widetilde{N}_0$, where $\widetilde{N}_0$ is the first integer (which exists and is finite) such that
$$
\forall j\geq \widetilde{N}_0, \quad \Vert \theta_j^T\Vert_{L^2(0,T)}^2\leq \frac{e^2\left(\pi L-\sin (\pi L)\right)}{128}e^{2T(j^2-1)}.
$$
As a result, $N_0$ is equal to $1$ if $T$ is large enough.
\end{theorem}

In other words, Theorem \ref{thm:stat1D} says that the sequence $(\omega^N)_{N\in\N^*}$ of optimal sets, whose existence is stated in Proposition \ref{truncTheo}, is stationary. The numerical procedure consists of computing these sets, and once it has become stationary, then we have found the optimal set $\omega^*$, solution of the shape optimization problem \eqref{def:odp1D}. 

A natural issue concerns the characterization of the minimal integer $N_0$ such that the sequence of optimal domains $(\omega^N)_{N\geq N_0}$ remains constant. Even if a partial answer is provided in Theorem \ref{thm:stat1D}, it is likely that the determination of $N_0$ is in general intricate. 

\begin{figure}[h!]
\begin{center}
\includegraphics[width=3.5cm]{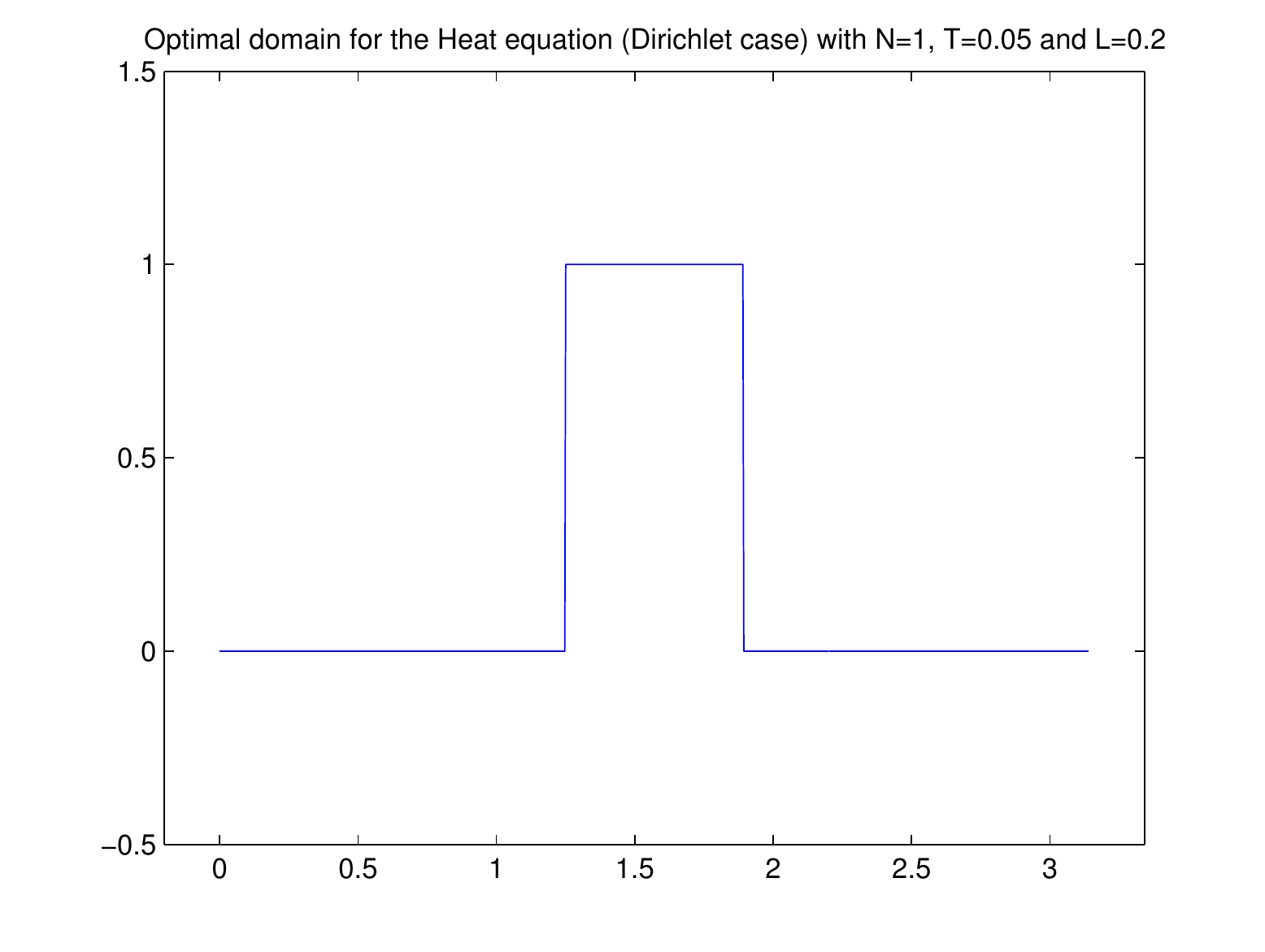}
\includegraphics[width=3.5cm]{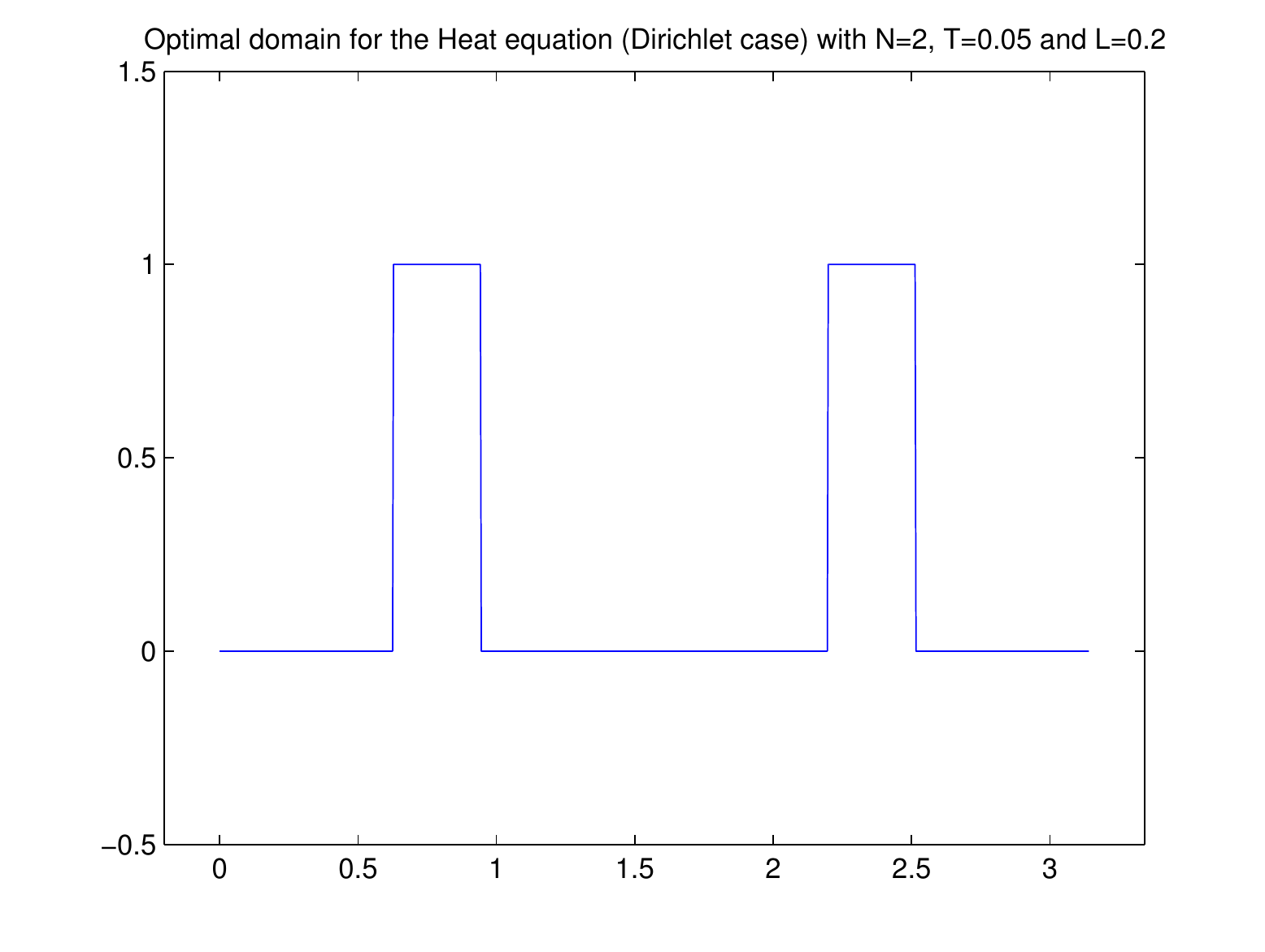}
\includegraphics[width=3.5cm]{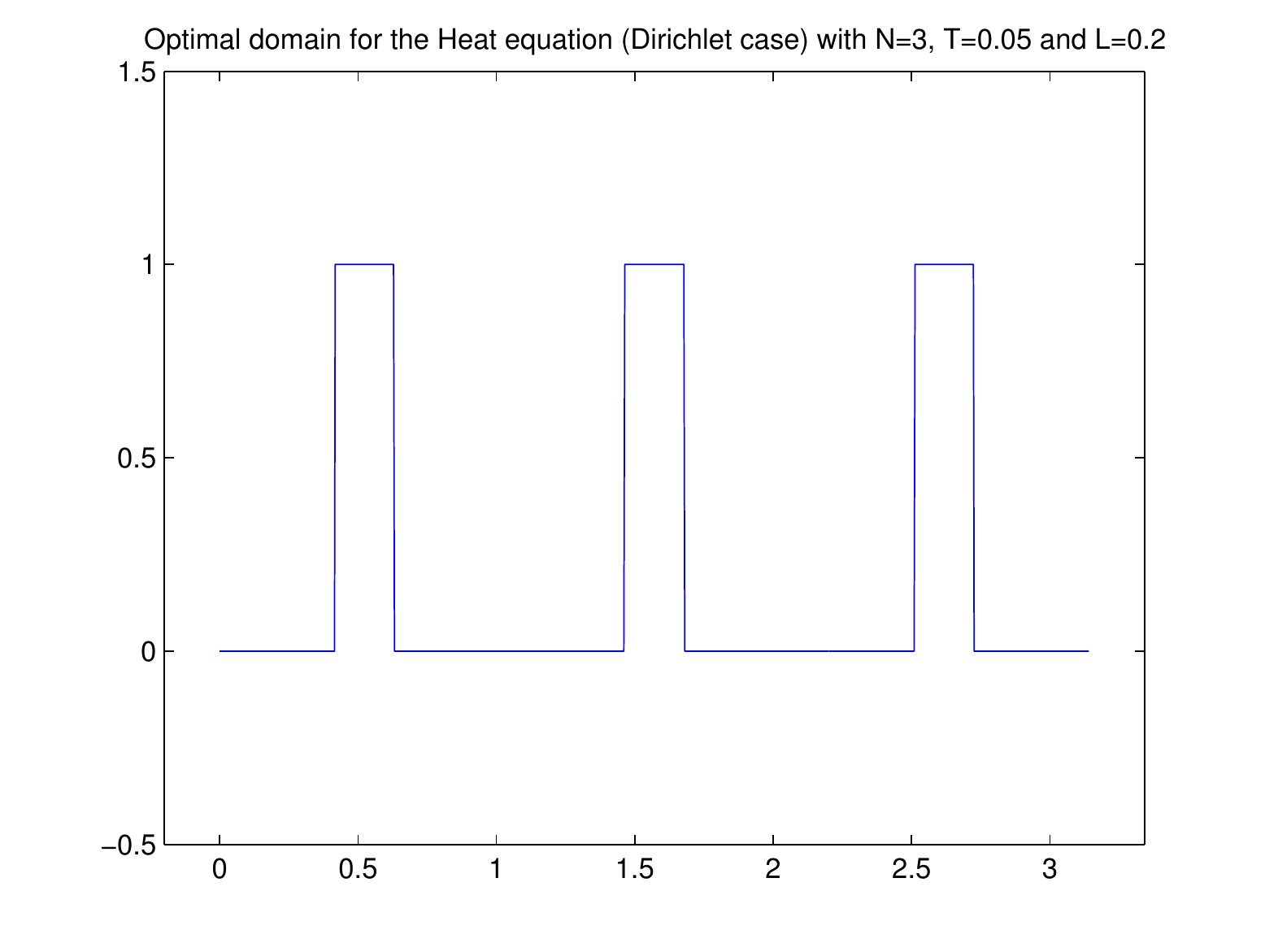}
\includegraphics[width=3.5cm]{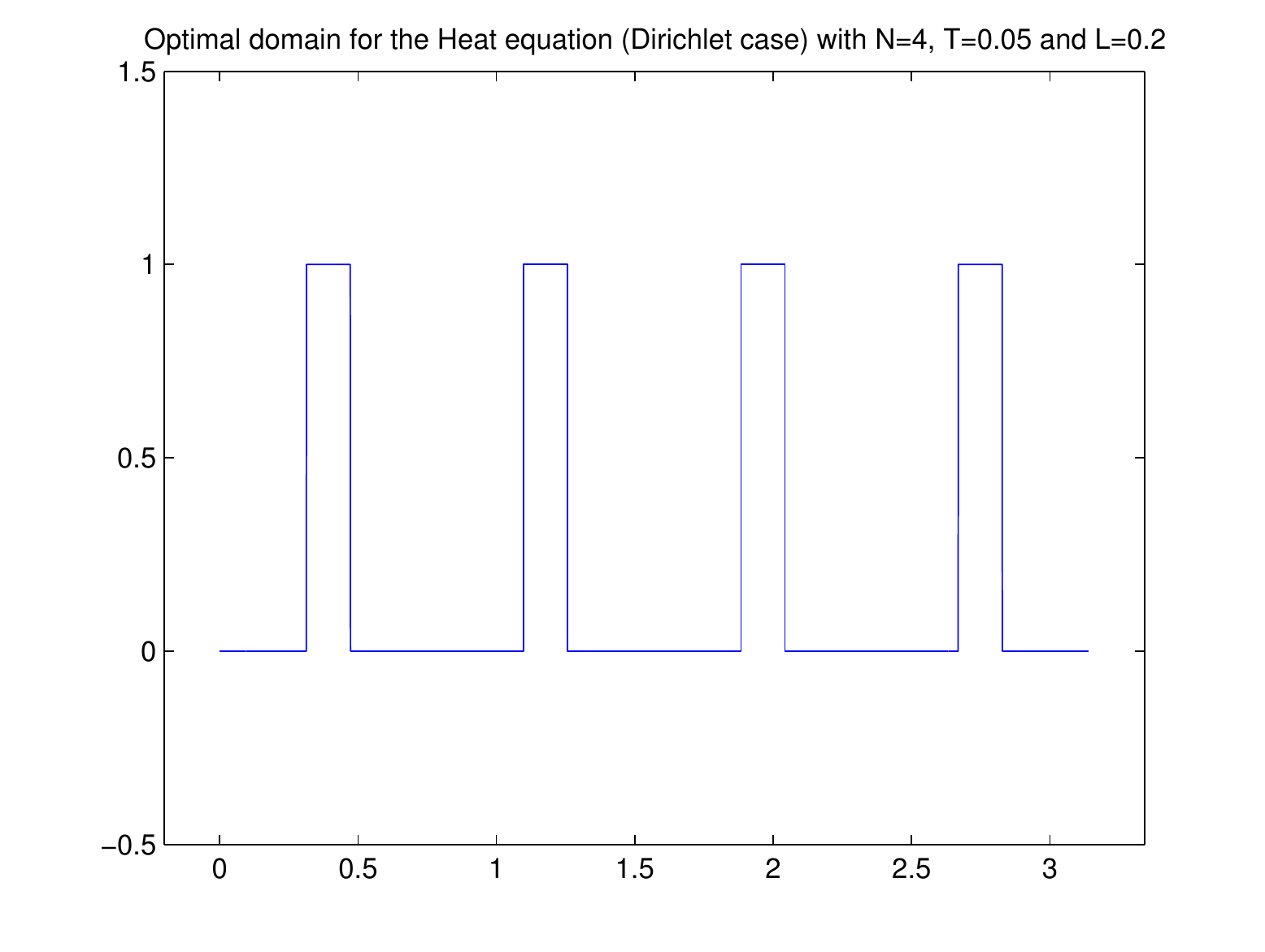}
\includegraphics[width=3.5cm]{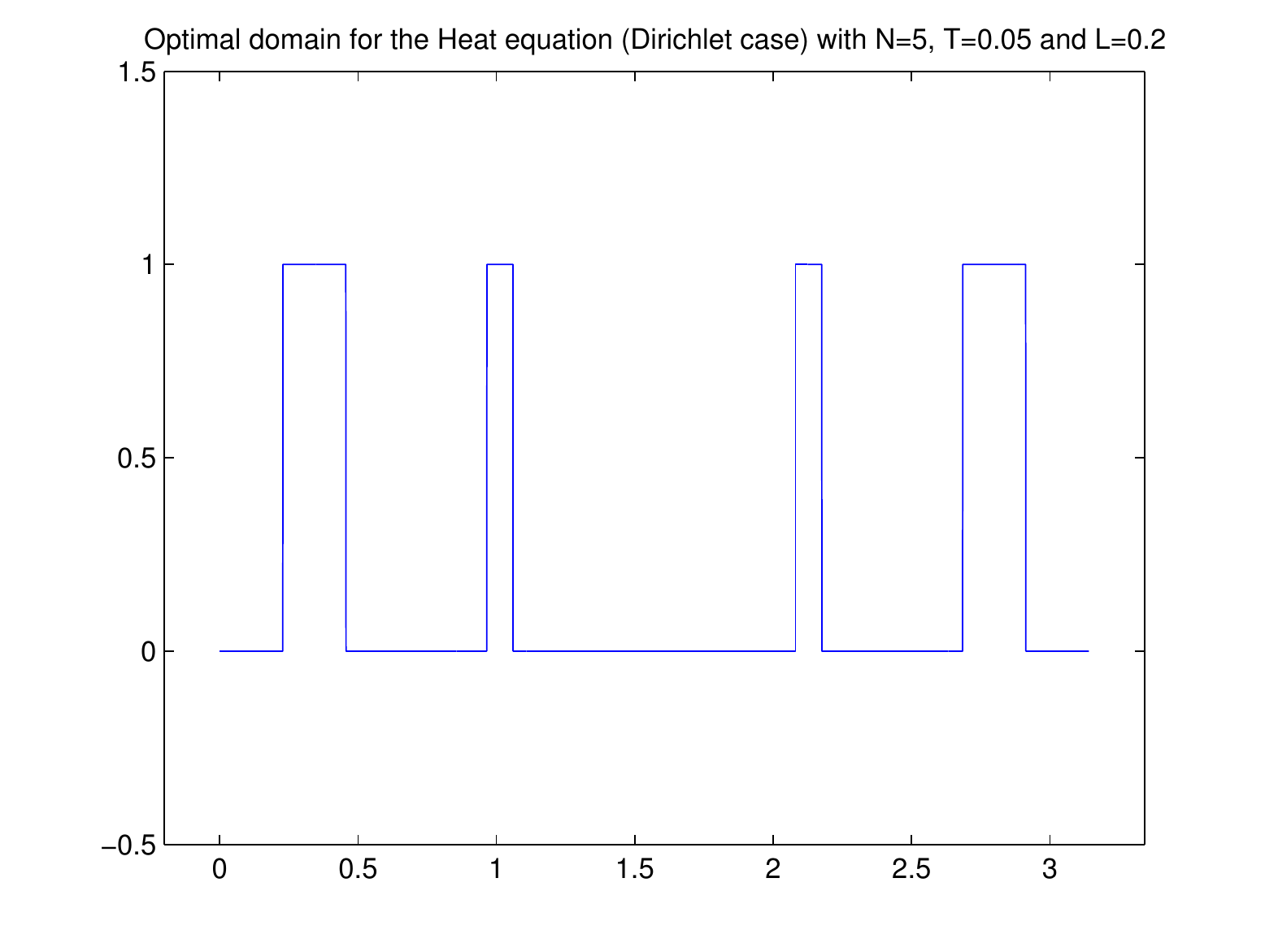}
\includegraphics[width=3.5cm]{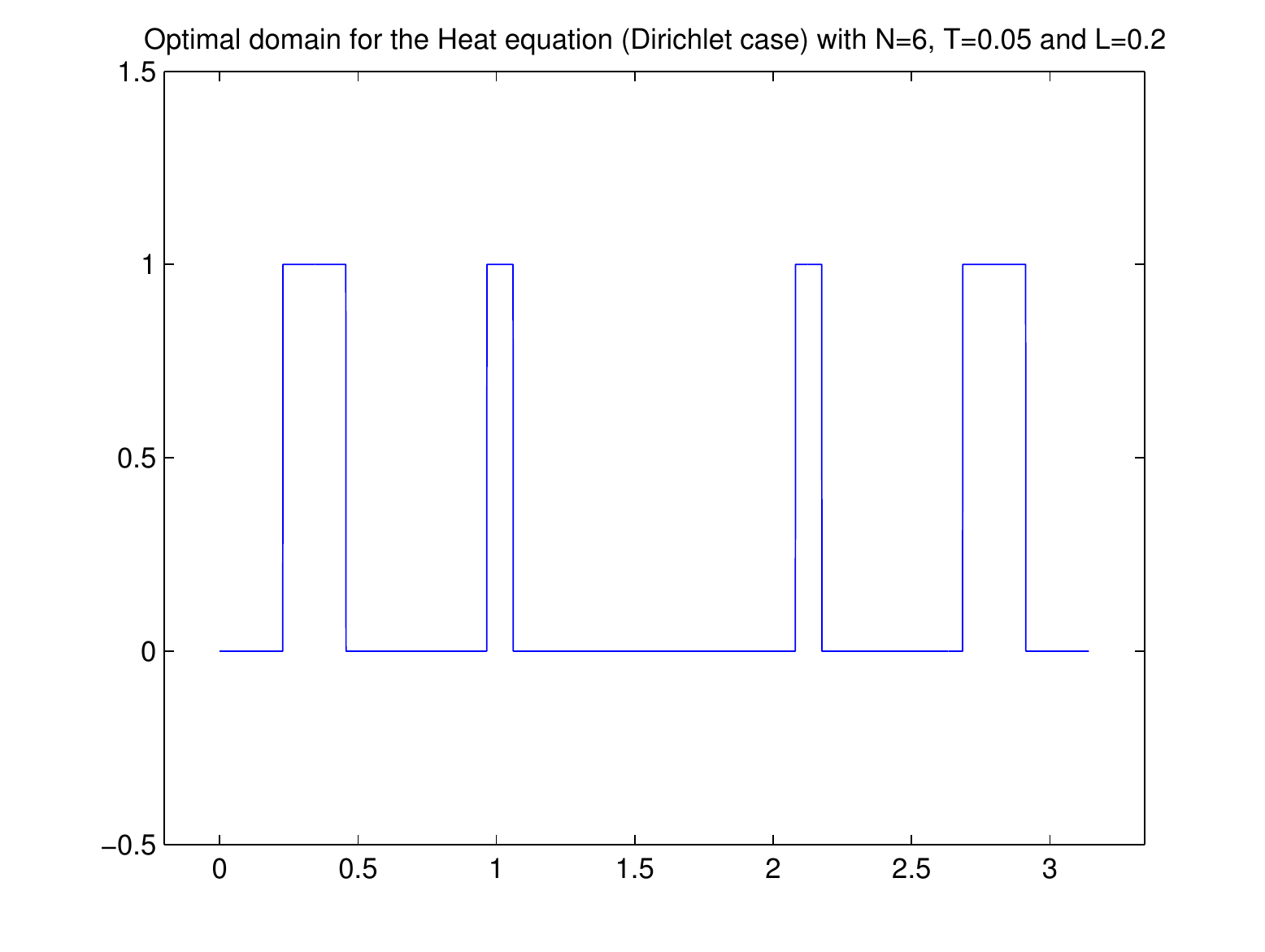}
\includegraphics[width=3.5cm]{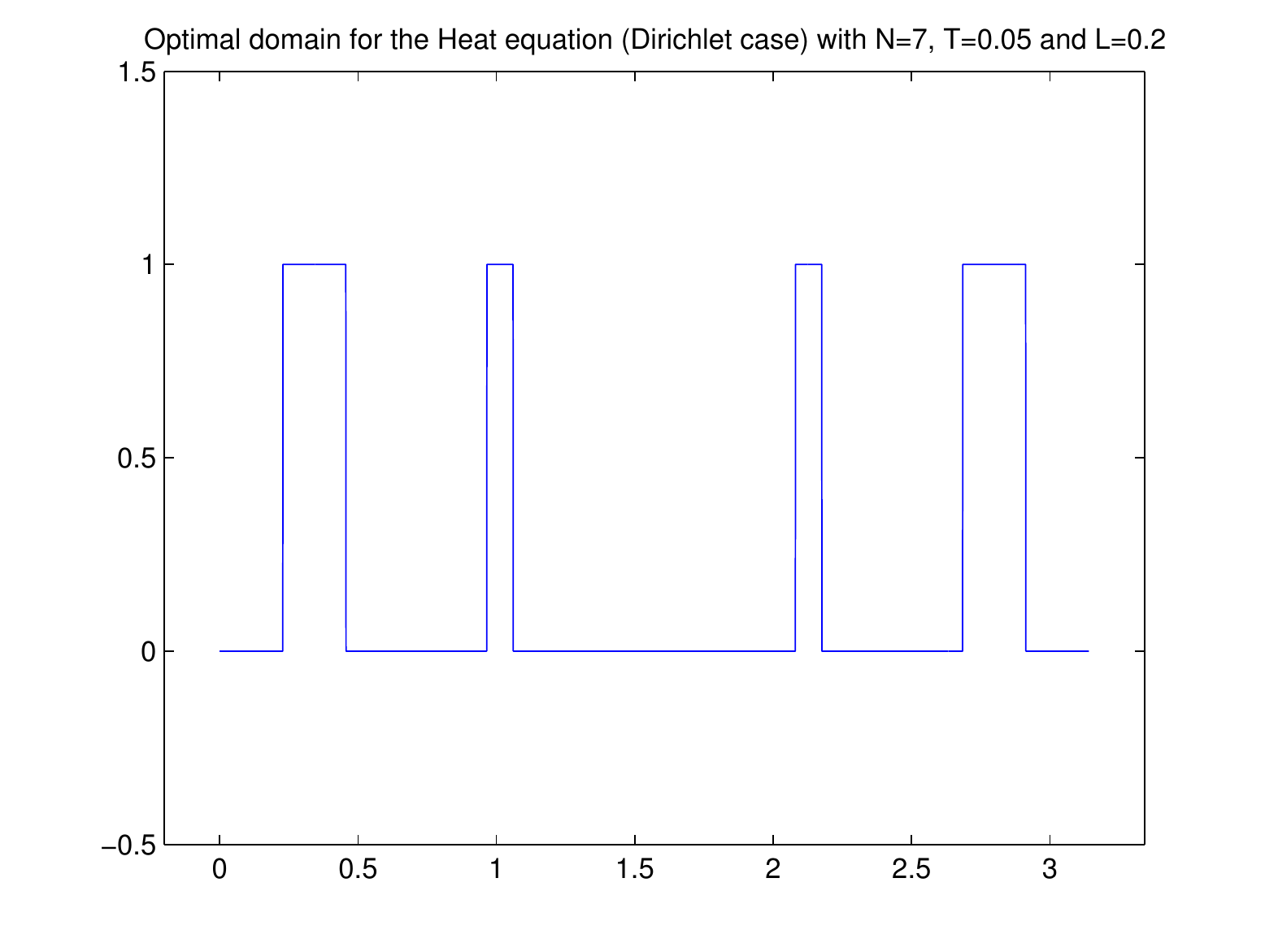}
\includegraphics[width=3.5cm]{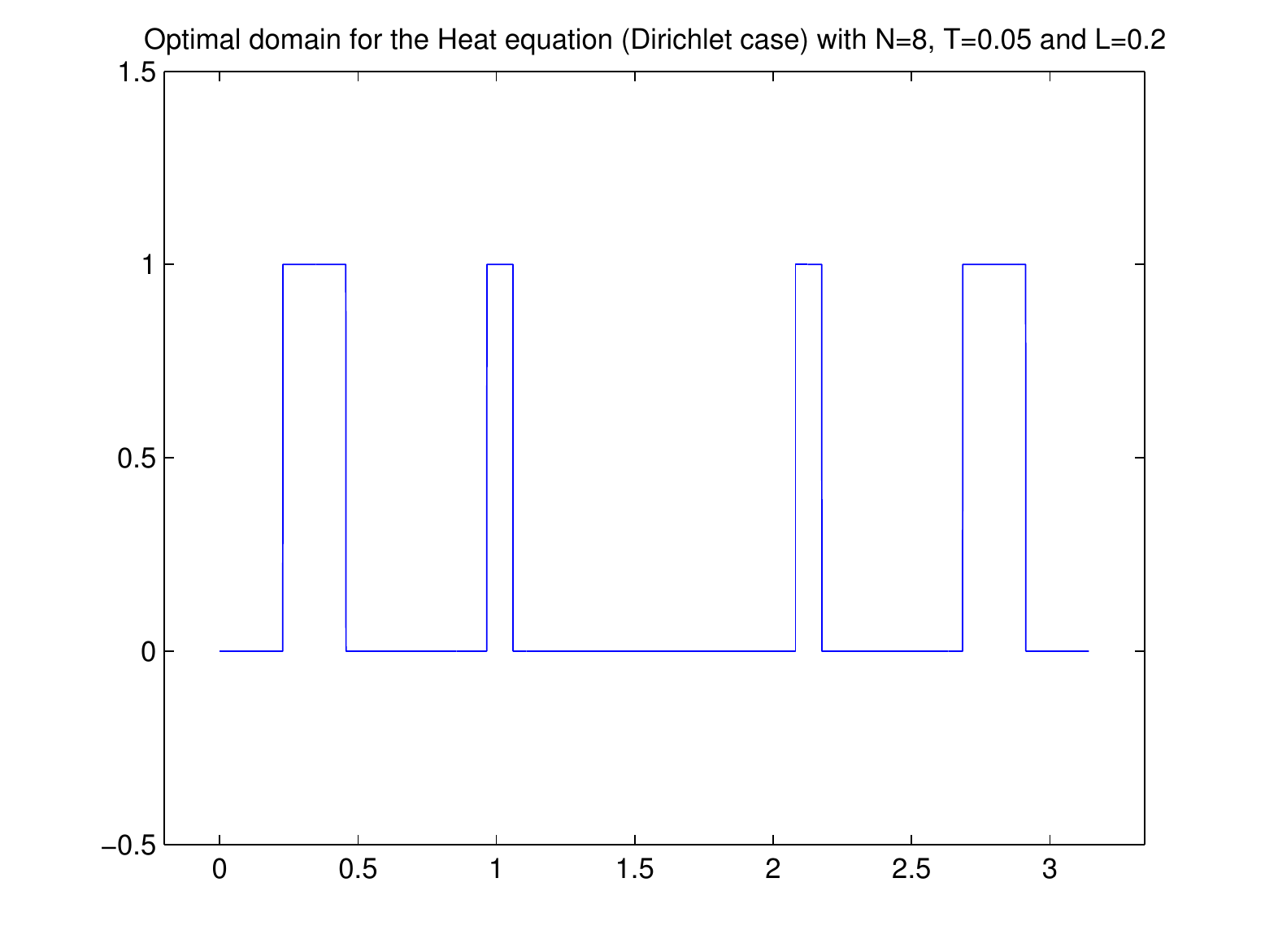}
\caption{${\Omega}=(0,\pi)$, $L=0.2$, $T=0.05$. 
From left to right, and top to down: optimal solution $\chi_{\omega^N}$  for $N=1,\ldots,8$.}\label{figpb2Dcont}
\end{center}
\end{figure}

As a numerical illustration of this computation procedure, we provide on Figure \ref{figpb2Dcont} several numerical simulations of the optimal control domain, solution of the truncated problem \eqref{pb_max_JN} in the 1D case, for the Dirichlet-Laplacian. We observe the expected stationarity property of the sequence of optimal domains $\omega^N$ from $N=5$ on.

\bigskip

In the forthcoming section devoted to providing the proofs of the results above, it will be required to consider a convexified version of the problem \eqref{def:odp1D}, which may fail to have some solutions because of the hard constraint\footnote{Indeed, equality constraints in $L^\infty$ are in general not preserved by the natural topologies such as the $L^\infty$ weak-star topology.} $\chi_\omega\in\mathcal{U}_L$ (which is a binary constraint almost everywhere). This is usually referred to as relaxation (see, e.g., \cite{BucurButtazzo}). Since the set $\mathcal{U}_L$ (defined by \eqref{defUL}) does not share nice compactness properties, we consider the convex closure of $\mathcal{U}_L$ for the weak star topology of $L^\infty$, which is
\begin{equation}\label{defALbar}
\overline{\mathcal{U}}_L = \left\{ a\in L^\infty(0,\pi;[0,1])\ \vert\ \int_{\Omega} a(x)\,dx=L\pi\right\}.
\end{equation}
Such a relaxation was used as well in \cite{munchHeat,PTZObs1,PTZobsND}.
Replacing $\chi_\omega\in\mathcal{U}_L$ with $a\in\overline{\mathcal{U}}_L$, we consider the relaxed (or convexified) formulation of the problem \eqref{def:odp1D} given by
\begin{equation}\label{defJa}
\sup_{a\in \overline{\mathcal{U}}_L} \mathcal{J}(a) ,
\end{equation}
where the functional $J$ is naturally extended to $\overline{\mathcal{U}}_L$ by
\begin{equation}\label{defJrelax}
\mathcal{J}(a)=\inf_{j\in \N^*}\frac{e^{2j^2T}}{\int_{0}^T\theta_j^T(t)^2\, dt} \int_{\Omega}a(x) \sin^2(jx)\, dx,
\end{equation}
for every $a\in \overline{\mathcal{U}}_L$. 
We consider as well a relaxed formulation of the truncated optimal design problem \eqref{pb_max_JN} by
\begin{equation}\label{defJaN}
\sup_{a\in \overline{\mathcal{U}}_L} \mathcal{J}_N(a) ,
\end{equation}
where the functional $\mathcal{J}_{N}$ is naturally extended to $\overline{\mathcal{U}}_L$ by
\begin{equation}\label{defJaNrelax}
\mathcal{J}_N(a)=\inf_{1\leq j\leq N}\frac{e^{2j^2T}}{\int_{0}^T\theta_j^T(t)^2\, dt} \int_{\Omega}a(x)\sin^2(jx)\, dx,
\end{equation}
for every $a\in \overline{\mathcal{U}}_L$. 

Being defined as the infimum of linear functions, continuous for the $L^\infty$ weak star topology, the functional $\mathcal{J}$ is upper semi continuous for the $L^\infty$ weak star topology. The set $\overline{\mathcal{U}}_L$ being compact for this topology, we then have the following result.

\begin{lemma}\label{propExistRelax}
For every $L\in (0,1)$, the relaxed problem \eqref{defJa} (respectively \eqref{defJaN}, for any $N\in\N^*$) has at least one solution $a^*\in\overline{\mathcal{U}}_L$ (respectively $a^{N}\in\overline{\mathcal{U}}_L$).
\end{lemma}

\subsection{Proof of Proposition \ref{truncTheo}}\label{sec_A2}
Considering the functions $a(\cdot)$ of $\overline{\mathcal{U}}_L$ as controls, and interpreting the problem \eqref{pb_max_JN} as an optimal control problem, leads to consider the control system
\begin{equation}\label{contsys1henrot}
\begin{split}
y'(x) &= a(x), \\
y_j'(x) & = \frac{e^{2j^2T}}{\int_{0}^T\theta_j^T(t)^2\, dt} a(x)\sin^2(jx),\quad j\in\{1,\ldots,N\}, \\
z'(x) & = 0, 
\end{split}
\end{equation}
for almost every $x\in[0,\pi]$, with initial conditions
\begin{equation}\label{initcond1henrot}
y(0)=0,\quad y_j(0)=0,\ j\in\{1,\ldots,N\}.
\end{equation}

The additional function $z$ above stands for the cost functional $\mathcal{J}_N(a)$ and will be defined with the help of inequality constraints below since it is written as the minimum of the quantities $y_j (\pi)$ over $ j\in\{1,\ldots,N\}$.

The relaxed problem \eqref{defJa} is then equivalent to the optimal control problem of determining a control $a\in\overline{\mathcal{U}}_L$ steering the control system \eqref{contsys1henrot} from the initial conditions \eqref{initcond1henrot} to the final condition
\begin{equation}\label{finalcondhenrot}
y(\pi)=L\pi,
\end{equation}
and maximizing the quantity $z (\pi)$ (or similarly $z(0)$, since $z$ in constant on $[0,\pi]$), with the additional final conditions
\begin{equation}\label{finalcondhenrot2}
z (\pi)\leq y_j (\pi),\quad \textnormal{for all }j\in\{1,\ldots,N\}.
\end{equation}
Indeed, this follows directly from the observation that the unique solution of 
$$
\max \{z\ \mid \ z\leq y_j (\pi), \ j\in\{1,\ldots,N\} \}
$$
is $z=\min_{1\leq j\leq N}y_j (\pi)$.

Therefore, $a^*$ is a solution of the optimal control problem above.
The existence of an optimal control is standard. According to the Pontryagin Maximum Principle (see \cite{Pontryagin}), if $a$ is optimal then there exist real numbers\footnote{Note that, since the dynamics of \eqref{contsys1henrot} do not depend on the state, it follows that the adjoint states of the Pontryagin Maximum Principle are constant.} $(p_y,p_1,\ldots,p_N)\in \R_{-}\times \R_{+}^{N}\backslash (0,\ldots,0)$, such that
\begin{equation}\label{maxcondhenrot}
a(x) = \left\{ \begin{split} 1 & \quad \textrm{if}\ \varphi^N(x)>0,\\ 0 &\quad  \textrm{if}\ \varphi^N(x)<0, \end{split} \right.
\end{equation}
for almost every $x\in[0,\pi]$, 
where the so-called switching function $\varphi^N$ is defined by
\begin{equation}\label{defvarphihenrot}
\varphi^N(x) = p_y+\sum_{j=1}^{N}\frac{e^{2j^2T}}{\int_{0}^T\theta_j^T(t)^2\, dt} p_j\sin^2(jx).
\end{equation}
Moreover, the control $a(\cdot)$ is nonsingular (see \cite{trelat}) since $\varphi^N$ is a finite trigonometric sum and thus cannot be constant on any subset of positive measure. In particular, this implies that the optimal control $a^N$ is the characteristic function of a measurable subset $\omega^N(L)$ of $[0,\pi]$ of measure $L\pi$.
Note that the minimum of $\varphi^N$ on $[0,\pi]$ is reached at $0$ and $\pi$, hence from \eqref{maxcondhenrot} the optimal set $\omega^N$ does not contain $0$ and $\pi$.


To prove uniqueness, according to the previous discussion where it is stated that every maximizer of $J$ over $\overline{\mathcal{U}}_{L}$  is the characteristic function of some subset of $[0,\pi]$,  assume that there exist two distinct minimizers $\chi_{\omega_1}$ and $\chi_{\omega_2}$ in $\mathcal{U}_{L}$. As a maximum of linear functionals, the functional $a\mapsto \mathcal{J}(a)$ is convex on $\overline{\mathcal{U}}_L$, and it follows that for every $t\in (0,1)$ the function $t\chi_{\omega_1}+(1-t)\chi_{\omega_2}$ is also a solution of the problem \eqref{defJaN}, which is in contradiction with the fact that any solution of this problem is extremal.

Finally, the fact that $\omega^N(L)$ has at most $N$ connected components follows from the facts that the elements of $\partial\omega^N(L)$ are the solutions of $\varphi^N(x)=0$ and that $\varphi^N$ can be written as
$$
\varphi^N(x)=p_y+\frac{1}{2}\sum_{j=1}^N\frac{e^{2j^2T}}{\int_{0}^T\theta_j^T(t)^2\, dt} p_j-\frac{1}{2}\sum_{j=1}^N \frac{e^{2j^2T}}{\int_{0}^T\theta_j^T(t)^2\, dt} p_{j}T_{2j}(\cos x),
$$
where $T_{2j}$ denotes the $2j$-th Chebychev polynomial of the first kind. The degree of the polynomial $\varphi^N(\arccos X)$ (in the argument $X$) is at most $2N$, whence the result.


\subsection{Proofs of Theorems \ref{thm:odp1D} and \ref{thm:stat1D}}\label{sec_thm:odp1D}
The main idea of this proof is close to the one of \cite[Theorem 1]{PTZparabND}. According to Lemma \ref{propExistRelax}, the relaxed optimal design problem \eqref{defJa} has at least one solution $a^*\in\overline{\mathcal{U}}_L$. 
We will prove simultaneously Theorems \ref{thm:odp1D} and \ref{thm:stat1D}, by showing that $a^*$ coincides with the solution $a^N$ of the truncated problem \eqref{pb_max_JN} for $N$ large enough.

First of all, as a consequence of \cite[Lemma 6]{PTZ_HUM1D}, we have
$
\int_{0}^\pi a^*(x)\sin^2(jx)\, dx\geq \frac{L\pi-\sin(L\pi)}{2},
$
and therefore,
\begin{equation}\label{metz1144}
\frac{e^{2j^2T}}{\int_{0}^T\theta_j^T(t)^2\, dt}\int_{0}^\pi a^*(x)\sin^2(jx)\, dx\geq \frac{e^{2j^2T}(L\pi-\sin(L\pi))}{2\int_{0}^T\theta_j^T(t)^2\, dt}
\end{equation}
for every $j\in \N^*$.
Besides, we have the following result on the growth of the biorthogonal sequence $(\theta_j^T)_{j\in\N^*}$, following from \cite[Theorem 3.2]{MicuZuazua}. 
\begin{lemma}\label{lemma:CTbiorth}
Let $T>0$. There exists $C_{T}>0$ such that
$$
C_{T} \int_{0}^T\theta_j^T(t)^2\, dt\leq e^{2\pi j},
$$
for every $j\in\N^*$.
\end{lemma}
It follows from this result that
$$
\frac{e^{2j^2T}}{\int_{0}^T\theta_j^T(t)^2\, dt}\geq C_{T}e^{2j^2T-2\pi j},
$$
for every $j\in \N^*$.

Combining these two facts, we infer that
\begin{equation}\label{eq:limproof}
\lim_{j\to +\infty}\frac{e^{2j^2T}}{\int_{0}^T\theta_j^T(t)^2\, dt}\int_{0}^\pi a^*(x)\sin^2(jx)\, dx=+\infty ,
\end{equation}
and moreover, there exists $N_0\in\N^*$ such that
\begin{equation}\label{propN0}
\inf_{j>N_0}\frac{e^{2j^2T}}{\int_{0}^T\theta_j^T(t)^2\, dt}\int_{\Omega} a^*(x) \sin^2(jx)\, dx > \frac{e^{2T}}{\int_{0}^T\theta_1^T(t)^2\, dt} .
\end{equation}
Since there holds in particular 
$$
\mathcal{J}_{N_0}(a^*)\leq \frac{e^{2j^2T}}{\int_{0}^T\theta_j^T(t)^2\, dt}\int_{\Omega} a^*(x)\sin^2(x)\, dx\leq \frac{e^{2j^2T}}{\int_{0}^T\theta_j^T(t)^2\, dt},
$$
we infer from \eqref{propN0} that
$$
\mathcal{J}(a^*)  = \min \left(\mathcal{J}_{N_0}(a^*),\inf_{j>N_0}\frac{e^{2j^2T}}{\int_{0}^T\theta_j^T(t)^2\, dt}\int_{\Omega} a^*(x) \sin^2(jx)\, dx\right) = \mathcal{J}_{N_0}(a^*) .
$$
Let us actually prove that  $\mathcal{J}(a^*) = \mathcal{J}_{N_0}(a^{N_0})$, where $a^{N_0}\in\mathcal{U}_L$ denotes the unique maximizer of $\mathcal{J}_{N_0}$, as stated in Lemma \ref{propExistRelax}. Since $a^{N_0}$ maximizes $\mathcal{J}_{N_0}$ over $\overline{\mathcal{U}}_L$, one has $\mathcal{J}(a^*) =\mathcal{J}_{N_0}(a^*) \leq \mathcal{J}_{N_0}(a^{N_0})$.
Let us argue by contradiction and assume that $\mathcal{J}_{N_0}(a^*) < \mathcal{J}_{N_0}(a^{N_0})$.
For every $t\in[0,1]$, we set $ a_t = a^* + t ( a^{N_0} - a^*)$.
Since $\mathcal{J}_{N_0}$ is concave (as an infimum of linear functionals), we get
$$
\mathcal{J}_{N_0}(a_t) \geq (1-t)\mathcal{J}_{N_0}(a^*) + t \mathcal{J}_{N_0}(a^{N_0}) > \mathcal{J}_{N_0}(a^*) = \mathcal{J}(a^*),
$$
for every $t\in(0,1]$, which means that
\begin{equation}\label{train17h39}
\inf_{1\leq j\leq N_0} \frac{e^{2j^2T}}{\int_{0}^T\theta_j^T(t)^2\, dt} \int_{\Omega} a_t(x) \sin^2(jx) \, dx > \inf_{1\leq j\leq N_0} \frac{e^{2j^2T}}{\int_{0}^T\theta_j^T(t)^2\, dt} \int_{\Omega} a^*(x) \sin^2(jx) \, dx \geq \mathcal{J}(a^*),
\end{equation}
for every $t\in(0,1]$. 
Besides, for every $\varepsilon>0$ there exists $t>0$ small enough such that
$$ 
\frac{e^{2j^2T}}{\int_{0}^T\theta_j^T(t)^2\, dt} \int_{\Omega} a_t(x) \sin^2(jx) \, dx  \geq (1-t)\frac{e^{2j^2T}}{\int_{0}^T\theta_j^T(t)^2\, dt} \int_{\Omega} a^*(x)\sin^2(jx) \, dx  \geq   \frac{e^{2T}}{\int_{0}^T\theta_1^T(t)^2\, dt} +\varepsilon,
$$
for every $j>N_0$. Therefore, 
\begin{equation}\label{train17h44}
\inf_{j> N_0} \frac{e^{2j^2T}}{\int_{0}^T\theta_j^T(t)^2\, dt}\int_{\Omega} a_t(x) \sin^2(jx) \, dx  > \frac{e^{2T}}{\int_{0}^T\theta_1^T(t)^2\, dt}.
\end{equation}
Since there holds in particular $\mathcal{J}_{N_0}(a_t)\leq \frac{e^{2T}}{\int_{0}^T\theta_1^T(t)^2\, dt}$, we infer from \eqref{train17h39} and \eqref{train17h44} that $\mathcal{J}(a_t) = \mathcal{J}_{N_0}(a_t) > \mathcal{J}(a^*)$, which contradicts the optimality of $a^*$.

Therefore $\mathcal{J}_{N_0}(a^*) =\mathcal{J}(a^*)= \mathcal{J}_{N_0}(a^{N_0})$, whence the result.

\paragraph{Estimate of the integer $N_0$.} It remains to provide an estimate for $N_0$.
We claim that any nonzero integer $\widetilde{N}_0$ such that the inequality
$$
\inf_{j>\widetilde{N}_0} \frac{e^{2j^2T}}{\int_{0}^T\theta_j^T(t)^2\, dt} \int_{\Omega} \chi_{\omega^N}(x) \sin(jx)^2\, dx >  \frac{e^{2T}}{\int_{0}^T\theta_1^T(t)^2\, dt},
$$
holds true satisfies $N_0\leq \widetilde{N}_0$ (in the sequel, we denote by $\widetilde{N}_0$ any integer such that the sequence $(\omega^N)_{N\geq \widetilde{N}_0}$ remains constant). 

To prove this claim, let us consider the simple case where $T\geq 1$. Notice that, in the next explanations, the lower bound on the time $T$ is not a restriction of our approach and can be chosen as small as desired with a slight adaptation of the following arguments. It is possible to perform more precise computations since in this case, we know at the same time several properties on the involved biorthogonal sequences $(\theta_j^T)_{j\in\N^*}$ as well as the useful spectral property: for all $j\in\N^*$, one has $\int_{0}^\pi \chi_{\omega^N}(x)\sin^2(jx)\, dx\geq \frac{L\pi-\sin(L\pi)}{2}$ according to \cite[Lemma 6]{PTZ_HUM1D}. As a consequence, following the proof of Theorem \ref{thm:odp1D} and by using in particular \eqref{metz1144}, $\widetilde{N}_0$ can be chosen to be any integer such that
$$
\forall j\geq \widetilde{N}_0, \quad \frac{e^{2j^2T}}{\Vert \theta_j^T\Vert_{L^2(0,T)}^2}\left(\frac{\pi L-\sin (\pi L)}{2}\right)\geq \frac{e^{2T}}{\Vert \theta_1^T\Vert_{L^2(0,T)}^2}.
$$
According to \cite[Theorem 3.2]{MicuZuazua}, there holds $\Vert \theta_1^T\Vert_{L^2(0,T)}^2\geq \frac{e^2}{64}$, and we infer that $\widetilde{N}_0$ can also be chosen such that 
$$
\forall j\geq \widetilde{N}_0, \quad \Vert \theta_j^T\Vert_{L^2(0,T)}^2\leq \frac{e^2\left(\pi L-\sin (\pi L)\right)}{128}e^{2T(j^2-1)}.
$$
It remains to provide an upper bound of the quantity $\theta_j^T$ for any $j\in\N^*$. To this aim, we will use that for a given $j\in \N^*$, the mapping $v:[0,T]\ni t\mapsto \frac{(-1)^j}{2j}e^{-j^2T}\theta_j^T(t)$ is the control of minimal $L^2(0,T)$-norm for the boundary control problem of steering the system 
\begin{equation}\label{eqphicont}
\begin{split}
& \partial_t \varphi(t,x)-\partial_{xx} \varphi (t,x)=0, \qquad (t,x)\in (0,T)\times(0,\pi),\\
& \varphi(t,0)=0, \quad \varphi (t,\pi)=v(t), \qquad\qquad t\in (0,T),
\end{split}
\end{equation}
with initial datum $\varphi(0,x)=\sin(jx)$ to zero in time $T$, as highlighted in \cite[Proposition 2.2]{MicuZuazua}. Consider the particular control function $w_j$ vanishing on the time interval $[0,T-1]$ and equal to the control constructed with the help of the Hilbert Uniqueness Method (HUM) for the control problem of steering the system \eqref{eqphicont} with initial datum $\varphi(0,x)=e^{-j^2(T-1)}\sin (jx)$ to zero in time $1$. The existence of $w_j$ is well-known and we refer for instance to \cite{TucsnakWeiss,zuazua}. More generally, the controllability property of \eqref{eqphicont} also implies the existence $C>0$ that does not depend on $j$ nor $T$ such that $\Vert w_j\Vert_{L^2(0,T)}\leq Ce^{-(T-1)j^2}$ for all $j\in\N^*$ since the sequence of functions $(x\mapsto\sin(jx))_{j\in\N^*}$ is uniformly bounded in $L^2(0,\pi)$. We thus infer that
$$
\Vert \theta_j^T\Vert_{L^2(0,T)}\leq 2Cje^{j^2}
$$
for every $j\in\N^*$, and therefore, it suffices to choose $\widetilde{N}_0$ such that
$$
j\geq \widetilde{N}_0\ \Rightarrow\ 2Cje^{j^2}\leq \frac{e^2\left(\pi L-\sin (\pi L)\right)}{128}e^{2T(j^2-1)}.
$$
This estimate shows in particular that $\widetilde{N}_0$ and thus $N_0$ are equal to $1$ if $T$ is large enough.

\section{Generalization to parabolic distributed parameter systems, and lumped control}\label{sec:Control}
In this section, we generalize the results obtained previously for the one-dimensional heat equation, to a large family of parabolic systems. In a second step, we consider an alternative way of acting on the system, by means of lumped controls. 

\subsection{Problem setting}\label{sec31}
Let $n\in\N^*$ be an integer, and let $\Omega$ be a bounded open connected subset of $\R^n$.
We consider the internally controlled parabolic distributed parameter system
\begin{equation}\label{heatEqcontrolled}
\partial_t y+A_{0} y=\chi_\omega u, \quad t\in (0,T),
\end{equation}
where $A_0:D(A_{0})\rightarrow L^2(\Omega,\C)$ is a densely defined operator that generates a strongly continuous semigroup on $L^2(\Omega,\C)$, $u\in L^2((0,T)\times \Omega,\C)$ is the control function, and $\omega\subset\Omega$ is a measurable subset standing for the control domain.

We assume that there exists an orthonormal basis $(\phi_j)_{j\in\N^*}$ of $L^2(\Omega,\C)$ consisting of eigenfunctions of $A_{0}$, associated with (complex) eigenvalues $(\lambda_j)_{j\in\N^*}$ such that $\Real(\lambda_1)\leq\cdots\leq \Real(\lambda_j)\leq\cdots$.

The one-dimensional heat equation investigated previously enters into this frame, but now the setting is much more general.
\medskip

The objective of this section is to give a precise sense to the question of optimizing the control domain $\omega$. As a first remark, let us note that, since the equation is parabolic and thus has smoothing properties, we focus on the exact null controllability problem, that is the problem of steering the system from any initial condition (in an appropriate functional space) to zero, within a time $T>0$. 

We use the moment method in order to derive a relevant model of optimal sensor shape and location with results valuable for almost every initial data.
This method provides a way of constructing a control achieving exact null controllability, for some given initial data $y^0\in L^2(\Omega)$. 
As explained below, this approach suffers however from restrictions related to the M\"untz-Sz\'asz theorem, and then cannot be applied to any parabolic system.

We address this control problem in the framework developed in \cite{fattorini2} (see also the survey \cite{Russell}) where the controllability problem is reduced to a moment problem which is solved explicitly with the help of a biorthogonal sequence to the family of exponential functions $\Lambda = (e^{-\lambda_{j}t })_{j\geq 1}$. 

Consider the control system \eqref{heatEqcontrolled} with the initial data 
\begin{equation}\label{y0}
y(0)=y^0=\sum_{j\in\N^*} a_j\phi_j \in L^2(\Omega).
\end{equation}
The moment method provides a control steering the parabolic system \eqref{heatEqcontrolled} to zero, as stated in the following result.

\begin{lemma}\label{lemma:momentCont}
We define formally the function $u$ by
\begin{equation}\label{defu}
u(t,x)=-\sum_{j\in\N^*} \frac{a_j e^{-\lambda_jT}}{\int_\omega|\phi_{j}(y)|^2 \, dy}\theta_j^T(T-t)\phi_j(x) ,
\end{equation}
for almost every $t\in (0,T)$ and every $x\in \Omega$. If this series defines a function of $L^2((0,T)\times\Omega)$, then this control is a solution of the problem of steering the system \eqref{heatEqcontrolled} from $y^0$ to $0$ in time $T$.
\end{lemma}
The proof of this lemma is done in Section \ref{secproof31}.

\begin{remark}\label{rem_Muntz}
Recall that such a biorthogonal sequence exists if and only if the family $\Lambda$ is \textit{minimal}, that is, every element $t\mapsto e^{-\lambda_jt}$ lies outside of the closure in $L^2(0,T)$ of the vector space spanned by all other elements $t\mapsto e^{-\lambda_kt}$, with $k\neq j$.
If this condition is fulfilled, then this biorthogonal sequence is uniquely determined if and only if the family $\Lambda$ is complete in $L^2(0,T)$.

It is well known, by the M\"untz-Sz\'asz theorem, that the family $\Lambda$ is complete in $L^2(0,T)$ (but not independent) if and only if 
$$\sum_{j\in\N^*}\frac{1}{\Real(\lambda_j)+\lambda}=+\infty,$$
for some real number $\lambda$ such that $\Real(\lambda_j)+\lambda>0$ for every $j\in\N^*$ (for instance, $\lambda=-\Real (\lambda_1)+1$ is suitable).
On the contrary, if this series is convergent then the closure of the span of $\Lambda$ is a proper subspace of $L^2(0,T)$, moreover $\Lambda$ is minimal and thus a biorthogonal sequence exists. 

Then, here, we are led to assume that the series is convergent, which is a quite strong restriction on the parabolic system under consideration.
\end{remark}

For every $y^0\in L^2(\Omega)$, we set $\Gamma_{\omega}(y^0)=\chi_\omega u$, where $u$ is the control defined by \eqref{defu}, steering the system \eqref{heatEqcontrolled} from $y^0$ to $0$ in time $T$. This defines an operator $\Gamma_{\omega}:L^2(\Omega)\rightarrow L^2((0,T)\times\Omega)$, called the \textit{moment control operator}, which is linear and continuous. Its norm is
$
\Vert \Gamma_{\omega}\Vert = \sup \{ \Vert \Gamma_{\omega} (y^0)\Vert_{L^2((0,T)\times\Omega)} \mid \Vert y^0\Vert_{L^2(\Omega)}=1 \}.
$

As in the previous section, we randomize the Fourier coefficients of a given $y^0\in D(A_{0})$, with $y^0=\sum_{j=1}^{+\infty}a_j\phi_j$, by defining $a_j^\nu=\beta^\nu_{j}a_j$ for every $j\in\N^*$, where $(\beta_{j}^\nu)_{j\in\N^*}$ is a sequence of independent real-valued random variables on a probability space $(\mathcal{X},\mathcal{A},\mathbb{P})$ having mean equal to $0$, variance equal to $1$, and a super exponential decay (for instance, independent Bernoulli random variables). Then we define
$$
\mathcal{K}(\chi_\omega) = \sup_{\Vert y^0\Vert_{L^2(\Omega)}=1} \mathbb{E}\left(\Vert \Gamma_{\omega}(y^0_{\nu})\Vert^2_{L^2((0,T)\times \Omega)}\right),
$$
where $y^0_\nu= \sum_{j=1}^{+\infty}\beta_j^\nu a_j \phi_j$, and $\mathbb{E}$ is the expectation over the space $\mathcal{X}$ with respect to the probability measure $\mathbb{P}$.

\begin{lemma}\label{lemmContrand}
There holds
$$
\mathcal{K}(\chi_\omega) = \left( \inf_{j\in\N^*}\gamma_{j}(T)\int_{\omega}|\phi_j(x)|^2\, dx \right)^{-1} ,
$$
where the coefficients $\gamma_{j}(T)$ are defined by
\begin{equation}\label{defgammajcontrol}
\gamma_{j}(T)=\frac{e^{2\Real(\lambda_{j})T}}{\int_{0}^T\theta_j^T(t)^2\, dt},
\end{equation}
for every $j\in\N^*$.
\end{lemma}

This lemma is proved in Section \ref{sec_proof_lem2}.
As discussed previously, we model the best actuator shape and placement problem as the problem of minimizing $\mathcal{K}$ over the set $\mathcal{U}_{L}$ defined by
\begin{equation}\label{defULnew}
\mathcal{U}_L=\{ \chi_\omega\in L^\infty(\Omega ,\{0,1\})\ \vert\ \omega\subset \Omega\ \textrm{measurable} ,\  \vert\omega\vert=L |\Omega|\}.
\end{equation}
According to Lemma \ref{lemmContrand}, the problem of optimal actuator placement is equivalent to the problem
\begin{equation}\label{optDesignContrand}
\boxed{
\sup_{\chi_{\omega}\in \mathcal{U}_{L}}\inf_{j\in\N^*}\gamma_{j}(T)\int_{\omega}|\phi_j(x)|^2\, dx,
}
\end{equation}
where the coefficients $\gamma_{j}(T)$ are defined by \eqref{defgammajcontrol}.
In what follows, we define
$$
\mathcal{J}(\chi_\omega)=\inf_{j\in\N^*}\gamma_{j}(T)\int_{\omega}|\phi_j(x)|^2\, dx,
$$ 
for every measurable subset $\omega\subset\Omega$. 

\subsection{Main result and examples}\label{sec32}
We consider the following assumptions.
\begin{itemize}
\item[$\Hun$] (\textit{Strong Conic Independence Property}) If there exist a subset $E$ of $\Omega$ of positive Lebesgue measure, an integer $N\in\N^*$, a $N$-tuple $(\alpha_{j})_{1\leq j\leq N}\in (\R_{+})^N$, and $C\geq 0$ such that $\sum_{j=1}^N \alpha_j \vert\phi_j(x)\vert^2 = C$ almost everywhere on $E$, then there must hold $C = 0$ and $\alpha_j = 0$ for every $j\in \{1,\cdots,N\}$.
\item[$\Hdeux$] For every $a\in L^\infty(\Omega;[0,1])$ such that $\int_{\Omega}a(x)\, dx = L|\Omega|$, one has
\begin{equation*}
\liminf_{j\rightarrow+\infty} \ \gamma_j(T) \int_\Omega a(x) \vert\phi_j(x)\vert^2\, dx > \gamma_1(T) ;
\end{equation*}
\item[$\Htrois$] The eigenfunctions $\phi_j$ are analytic in $\Omega$.
\end{itemize}
These assumptions have been considered as well in \cite{PTZparabND} and are commented in that reference. For instance, they are satisfied for $A_0=(-\triangle)^\alpha$ with $\alpha>\frac{1}{2}$ and $\triangle$ is the Dirichlet-Laplacian on a piecewise $C^1$ domain $\Omega$ (see \cite[Section 2.4]{PTZparabND}).

The problem \eqref{optDesignContrand} is similar to the optimal design problem \eqref{def:odp1D}, except that now the weights $\gamma_j(T)$ are defined by \eqref{defgammajcontrol}.
It appears then important to estimate the asymptotics of $\gamma_j(T)$ as $j$ tends to $+\infty$. But this has been done in \cite{Avdonin,fattorini2,Hansen,MicuZuazua}. Those estimates will impose further restrictions on the problem under consideration.

For every $N\in\N^*$, we define the truncated criterion
\begin{equation*}
\mathcal{J}_N(\chi_\omega)=\inf_{1\leq j\leq N}\gamma_{j}(T)\int_{\omega}|\phi_j(x)|^2\, dx,
\end{equation*}
for every measurable subset $\omega\subset\Omega$. 
We have the following result.

\begin{proposition}\label{prop:trunc}
Let $N\in\N^*$. Under $\Hun$, the problem
\begin{equation}\label{pb:truncated}
\sup_{\chi_{\omega}\in\mathcal{U}_{L}}\mathcal{J}_N(\chi_\omega)
\end{equation}
has a unique solution $\chi_{\omega^N}$ in $\mathcal{U}_{L}$. Moreover, under the additional assumption $\Htrois$, $\omega^N$ is an open semi-analytic\footnote{A subset $\omega$ of a real analytic finite dimensional manifold $M$ is said to be semi-analytic if it can be written in terms of equalities and inequalities of analytic functions.
We recall that such semi-analytic subsets are stratifiable in the sense of Whitney (see \cite{Hardt,Hironaka}), and enjoy local finitetess properties, such that: local finite perimeter, local finite number of connected components, etc.} set.
\end{proposition}

This proposition is proved in Section \ref{secproof:prop:trunc}.
The main result is then the following theorem, proved in Section \ref{sec:proofcormaintheo}.

\begin{theorem}\label{cormaintheo}
Assume that there exist $m_1>0$, $m_2\in(0,2T)$, and a sequence $(\theta_j^T)_{j\in\N^*}$ biorthogonal to the family $\Lambda=(t\mapsto e^{-\lambda_{j}t})_{j\geq 1}$, such that
\begin{equation}\label{cond_vplambdaj}
\Vert \theta_j^T\Vert_{L^2(0,T)}^2 \leq m_1 e^{m_2 \Real(\lambda_j)},
\end{equation} 
for every $j\in\N^*$.
Then, under $\Hun$ and $\Hdeux$, the problem \eqref{optDesignContrand} has a unique solution $\chi_{\omega^*}\in\mathcal{U}_L$. Moreover there exists $N_{0}\in \N^*$ such that $\omega^*=\omega^N$, for every $N\geq N_0$. In particular, if $\Htrois$ is moreover satisfied, then $\omega^*$ is an open semi-analytic subset of $\Omega$, and thus, it has a finite number of connected components.
\end{theorem}

The same considerations as in Section \ref{sec:mainResIntro} on the algorithmic computation procedure still hold in this general framework. 

\bigskip

To finish, we provide hereafter some classes of examples for which the existence of a biorthogonal sequence satisfying \eqref{cond_vplambdaj} is known.

\begin{itemize}
\item Assume that there exist $\delta>0$, $\beta>1$, $\varepsilon>0$, $A\geq 0$ and $B\geq\delta$ such that
\begin{equation}\label{cond_vpj}
\vert \lambda_j-\lambda_k\vert \geq \delta \vert j^\beta-k^\beta \vert\quad\textrm{and}\quad
\varepsilon(A+Bj^\beta)\leq\vert\lambda_j\vert<A+Bj^\beta,
\end{equation}
for all $(j,k)\in(\N^*)^2$, where the elements of the sequence $(\lambda_k)_{k\in\N^*}$ are assumed to lie in $\{\lambda \in \C\mid |\arg \lambda |\leq \theta\}$ for some given $\theta\in (0,\pi/2)$. As argued in Remark \ref{rem_Muntz}, under the condition \eqref{cond_vpj} there exists a sequence $(\theta_j^T)_{j\in\N^*}$ biorthogonal to $\Lambda$, and it is proved in \cite{Hansen} that there exist two positive constants $\tilde A$ and  $\tilde B$ such that
$$
\Vert \theta_j^T\Vert_{L^2(0,T)}^2\leq \tilde B e^{\tilde A j},
$$
and since 
$$
\Real (\lambda_j)\geq |\lambda_j|\cos \theta\geq \varepsilon(A+Bj^\beta)\cos\theta
$$ 
for every $j\in\N^*$, we infer the existence of $m_1$ and $m_2$ such that the estimate \eqref{cond_vplambdaj} holds. We also refer to \cite[Theorem 3.2]{MicuZuazua} for an elementary proof of \eqref{cond_vplambdaj} for the eigenvalues the one-dimensional Dirichlet Laplacian operator.

For example, assume that $A_{0}=(-\triangle)^\alpha$ is a positive power of the one-dimensional Dirichlet-Laplacian on $\Omega=(0,\pi)$; then \eqref{cond_vpj} is satisfied if and only if $\alpha>1/2$.

In \cite{Hansen} other examples are provided where \eqref{cond_vpj} is satisfied, such as the damped Euler-Bernoulli plate in dimension two.

\item Assume that $(\lambda_n)_{n\in\N^*}$ is a sequence of positive real numbers and that there exist $K>0$, $\alpha>0$ and $\beta>1$ such that
$$
\lambda_{n}=K(n+\alpha)^\beta+\operatorname{o}(n^{\beta-1}),
$$
as $n$ tends to $+\infty$. 
It is proved in  \cite[Formula (3.25)]{fattorini2} that there exists two constants $\tilde A$ and $\tilde B$ such that 
$$
\Vert \theta_j^T\Vert_{L^2(0,T)}^2\leq \tilde A e^{\tilde B \lambda_j^{1/\beta }}
$$
for every $j\in \N^*$ and the estimate \eqref{cond_vplambdaj} then holds true. Note that the authors of \cite{fattorini2} use it to derive exact controllability results for a Sturm-Liouville one-dimensional equation.
We also mention the article \cite{ammar} where the authors extend the above approach and estimate to the framework of systems of one-dimensional parabolic equations, in view of establishing exact boundary controllability properties.

\item Assume that $A_0$ is the Dirichlet-Laplacian on the unit ball $\Omega=\{x\in \R^n\mid \Vert x\Vert < 1\}$, with $n$ arbitrary.
Using a refined study of the sequences of eigenfunctions and eigenvalues, it is proved in \cite[Section 6, (6.27)]{fattorini1} that \eqref{cond_vplambdaj} holds true with a constant $m_{2}$ not depending on $T$, and the authors use it to investigate boundary controllability issues for the heat equation in $\Omega$. Then Theorem \ref{cormaintheo} can be applied, provided that $T$ is large enough (since it is required that $m_2\in (0,2T)$).
\end{itemize}

\subsection{Optimal lumped controls}\label{sec:OLC}
In this section, we investigate a variant of the previously studied optimal design problem, based on another kind of controls referred to in the literature as the \textit{lumped controls} (see \cite[Chapter 4]{Russell} or \cite[Chapter 1.4]{Khapalov}). This wording designates tensorized controls that are the product of separated variables functions in time and space, the space profile of the control term being given. Then one only acts on the system by means of tuning the time-intensity of the control.

Let $\Omega$ be an open connected subset of $\R^n$ and $A_0:D(A_{0})\rightarrow L^2(\Omega,\C)$ be a densely defined operator that generates a strongly continuous semigroup on $L^2(\Omega,\C)$. We adopt the same framework as in Section \ref{sec31}, assuming the existence of an orthonormal basis $(\phi_j)_{j\in\N^*}$ of $L^2(\Omega,\C)$ consisting of eigenfunctions of $A_{0}$, associated with (complex) eigenvalues $(\lambda_j)_{j\in\N^*}$ such that $\Real(\lambda_1)\leq\cdots\leq \Real(\lambda_j)\leq\cdots$.

Consider the internally controlled parabolic system
\begin{equation}\label{heatEqcontrolled_lumped}
\partial_t y(t,x)+A_{0} y(t,x)+g(x)u(t)=0, \quad (t,x)\in (0,T)\times \Omega,
\end{equation}
with Dirichlet boundary conditions, where $g\in L^2(\Omega,\C)$ is the control profile and $u\in L^2(0,T)$ is the control function. 
The controlled system \eqref{heatEqcontrolled_lumped} is a particular version of \eqref{heatEqcontrolled_intro}.

In some sense, the function $g$ plays the role of $\chi_\omega$ in \eqref{heatEqcontrolled}, but here, the control function $u$ depends only on $t$. The function $g$ is usually fixed and the control is $u$. Here, we propose to optimize the control profile $g$.

Performing the same analysis as in Section \ref{sec:control2001} and using the same notations, one proves easily that every initial datum $y^0=\sum_{j=1}^{+\infty}a_{j}\phi_{j}\in L^2(\Omega)$ can be steered to zero in time $T$ with the control $u\in L^2(0,T)$ given by
$$
u(t)=-\sum_{j=1}^{+\infty} \frac{a_j e^{-\lambda_{j}T}}{\int_\Omega g(y)\overline{\phi}_j(y) \, dy}\theta_j^T(T-t),
$$
provided that the Fourier coefficients $\int_\Omega g(y)\phi_{j}(y) \, dy$ of $g$ do not vanish.

As previously, we define the \textit{moment control operator} $\tilde \Gamma_{g}:L^2(\Omega)\rightarrow L^2((0,T)\times \Omega)$ by $\tilde \Gamma_{g}(y^0)=f$, with $f(t,x)=g(x)u(t)$. Its norm is given by
$$
\Vert \tilde \Gamma_{g}\Vert =\sup_{\Vert y^0\Vert_{L^2(\Omega)}=1}\Vert \tilde \Gamma_{g}(y^0)\Vert _{L^2((0,T)\times\Omega )}=\Vert g\Vert_{L^2(\Omega)}\sup_{\Vert y^0\Vert_{L^2(\Omega)}=1}\Vert u\Vert_{L^2(0,T)}
$$
Following the framework developed in Sections \ref{sec:control2001} and \ref{sec31} leads to define a randomized criterion by defining $a_j^\nu=\beta^\nu_{j}a_j$ for every $j\in\N^*$. 
Then we define
$$
\tilde{\mathcal{K}}_g(\chi_\omega) =\sup_{\Vert y^0\Vert_{L^2(\Omega)}=1}\mathbb{E}(\Vert \tilde \Gamma_{g}(y^0_\nu)\Vert_{L^2((0,T)\times \Omega)}^2 ),
$$
where $y^0_{\nu}$ denotes the function of $L^2(\Omega)$ whose Fourier coefficients are the $a_{j}^\nu$ defined above. 

\begin{lemma}
There holds
$$
\tilde{\mathcal{K}}_g(\chi_\omega) = \sup_{j\in\N^*}\frac{e^{-2\Real(\lambda_{j})T}\int_{0}^T\theta_j^T(t)^2\, dt}{\left|\int_{\Omega} g(x)\overline{\phi}_{j} (x)\, dx\right|^2}\int_{\Omega} |g(x)|^2\, dx.
$$
\end{lemma}

The proof is similar to the proofs of Lemmas \ref{lemm:contHeat} and \ref{lemma:momentCont}, and thus is skipped.

We model the ``best design of lumped controller" as the problem of minimizing $\tilde{\mathcal{K}}_g(\chi_\omega)$ over the set of all possible profiles $g\in L^2(\Omega)$. The functional $g\mapsto \tilde{\mathcal{K}}_g(\chi_\omega)$ being homogeneous according to the previous lemma, the problem of optimal lumped control placement is then equivalent to the problem
\begin{equation}\label{pboptGamma1010}
\boxed{
\sup_{\Vert g\Vert_{L^2(\Omega)}=1}\ \inf_{j\in\N^*}\gamma_{j}(T)\left|\int_{\Omega} g(x)\overline{\phi}_{j} (x)\, dx\right|^2,
}
\end{equation}
with
$$
\gamma_{j}(T)=\frac{e^{2\Real(\lambda_{j})T}}{\int_{0}^T\theta_j^T(t)^2\, dt} ,
$$
for every $j\in \N^*$. Let us now solve this optimal design problem.

\begin{theorem}\label{thm:lumped}
We assume that
\begin{equation}\label{assump:specGammaj}
\sum_{j=1}^{+\infty}\frac{1}{\gamma_{j}(T)}<+\infty .
\end{equation}
Then, the problem \eqref{pboptGamma1010} has at least one solution, and we have
$$
\sup_{\Vert g\Vert_{L^2(\Omega )}=1}\ \inf_{j\in\N^*}\gamma_{j}(T)\left|\int_{\Omega} g(x)\overline{\phi}_{j} (x)\, dx\right|^2=\left(\sum_{j=1}^{+\infty}\frac{1}{\gamma_{j}(T)}\right)^{-1}.
$$
Moreover, the set of solutions consists of all functions $g$ in $L^2(\Omega)$ that can be expanded as
$$
g=\sum_{j=1}^{+\infty}g_{j}\phi_{j}\qquad\textrm{with}\qquad |g_{j}|^2=\left(\sum_{j=1}^{+\infty}\frac{1}{\gamma_{j}(T)}\right)^{-1}\frac{1}{\gamma_{j}(T)}.
$$
\end{theorem}

This theorem is proved in Section \ref{secproof:lumpedthm}.

\begin{remark}
Consider the case where $A_{0}=-\partial_{xx}$ is defined on $H^2(0,\pi)\cap H^1_{0}(0,\pi)$ (one-dimensional Dirichlet-Laplacian). Then $\phi_{j}(x)=\sqrt{\frac{2}{\pi}}\sin (jx)$ and $\lambda_{j}=j^2$ for every $j\in\N^*$. 

Denote by $g$ any solution of the problem \eqref{pboptGamma1010}. 
According to \cite[Theorem 3.2]{MicuZuazua}, there exists a positive constant $C_{T}$ such that
$$
\frac{e^{2j^2T}}{\int_{0}^T\theta_j^T(t)^2\, dt}\geq C_{T}e^{2j^2T-2\pi j} ,
$$
for every $j\in\N^*$. According to Theorem \ref{thm:lumped}, it follows that the Fourier coefficients $g_{j}$ decrease exponentially with respect to $j$, and as a consequence, the optimal functions $g$ are analytic (see e.g. \cite[Chapter 11, \S 63]{akhiezer}).
\end{remark}

\begin{remark}
It might seem natural and of physical interest to investigate what happens if we restrict our search of the control profile $g$ to a set of characteristic functions of a measurable subset $\omega$,  with the measure of $\omega$ possibly fixed. Doing this, we get a kind of instability: indeed, assuming that $\omega$ is the finite union of rational intervals (in other words, intervals whose extremities are rational multiples of $\pi$), one can easily check that
$$
\inf_{j\in\N^*}\gamma_{j}(T)\left(\int_{0}^\pi\chi_{\omega}(x)\sin (jx)\, dx\right)^2=0.
$$
Therefore, this problem appears to be ill-posed in some sense, and is probably not so much relevant with respect to practical issues.
\end{remark}

\subsection{Conclusion}\label{sec:conclpb}

To conclude, let us provide several further comments and open problems.

\paragraph{Generalization to other methods of control and higher dimensions.}
As underlined in the previous sections, the use of controllers obtained by the moment method reduces mainly the perimeter of our study to one-dimensional operators. 

In view of generalizing our approach to other control operators, let us use the framework described in Section \ref{sec31}, considering the controlled system
\begin{equation}\label{syst:concl}
\partial_t y+A_{0} y=\chi_\omega u, \quad t\in (0,T),
\end{equation}
where $y(0,\cdot)=y^0\in L^2(\Omega)$ and $u$ is a control steering this system from $y^0$ to 0 in time $T$, whenever it is possible. Let us assume that for every $T>0$ and every Lebesgue measurable subset $\omega$ of positive measure, the system \eqref{syst:concl} is null-controllable in time $T$. In this case, let us write $\Gamma_\omega=\chi_\omega u$.

For instance, the Hilbert Uniqueness Method (see \cite{lions2,lions}) is a well-known method used to design a null control for \eqref{heatEqcontrolled_intro}-\eqref{dir_cond}, with the additional property that this control has a minimal $L^2$ norm over all possible null controls. The null-controllability property in time $T$ of this system is equivalent to an observability property on the pair $(\omega,T)$. Note that in the case where $A_0$ is the Dirichlet-Laplacian operator $-\triangle$, it has been showed in \cite{AEWZ} that the observability inequality holds true for every $T>0$ and every Lebesgue measurable subset $\omega$ of positive measure. 

Following the approach described in Section \ref{sec31}, we define 
$$
\mathcal{K}(\chi_\omega) = \sup_{\Vert y^0\Vert_{L^2(\Omega)}=1} \mathbb{E}\left(\Vert \Gamma_{\omega}(y^0_{\nu})\Vert^2_{L^2((0,T)\times \Omega)}\right),
$$
where $y^0_\nu= \sum_{j=1}^{+\infty}\beta_j^\nu a_j \phi_j$, and $\mathbb{E}$ is the expectation over the space $\mathcal{X}$ with respect to the probability measure $\mathbb{P}$. Similar computations as those of Section \ref{sec_proof_lem2} enable to show that 
$$
\mathcal{K}(\chi_\omega)= \sup_{j\in\N^*} \big\Vert \Gamma_{\omega}(\phi_j) \big\Vert_{L^2((0,T)\times \Omega)}^2
$$

As discussed previously, we model the best actuator shape and placement problem as the problem of minimizing $\mathcal{K}$ over the set $\mathcal{U}_{L}$. Analyzing this optimal design problem does not seem easy since it requires to know fine regularity properties of each control function $u_j$ defined by $\Gamma_\omega(\phi_j)=\chi_\omega u_j$.

\paragraph{Analysis of the full control operator.} One of the main issues that remains to be developed is whether one can attack the problem of the optimal design of the controllers and actuators without the diagonalization procedure by randomization. The issue is then much harder to handle, as it occurs at the level of the observability problem. Note also that in that case, because of possible interactions of all modes, it is unclear how complex the optimal sets are.

Actually, if one defines the \textit{Gramian} operator $G_T$ as the infinite dimensional symmetric nonnegative matrix whose coefficient at row $j$ and column $k$ is given by $\int_0^T e^{(j^2+k^2)t}\, dt \int_\omega \sin(jx)\sin(kx)\, dx$, the operator norm $\Vert \Gamma_\omega\Vert$ is the inverse of the smallest eigenvalue of $G_T$. 
The randomization procedure consists in dropping the non-diagonal terms in $G_T$, by considering the inverse of smallest eigenvalue of $\operatorname{diag}(G_T)$.

Concerning the particular case of controllers given by the Hilbert Uniqueness Method (see \cite{lions2,lions}) to design a null control for \eqref{heatEqcontrolled_intro}-\eqref{dir_cond}, minimizing the control efforts in a deterministic way is actually equivalent to maximizing
$$
C_T(\chi_\omega)=\inf\left\{ \frac{ \int_0^T\int_\omega \vert y(t,x)\vert^2\,dx \, dt }{\Vert y(T,\cdot)\Vert^2_{L^2(0,\pi)}} \ \big\vert\  y(0,\cdot)\in L^2(0,\pi) \setminus\{0\} \right\} ,
$$
which is the largest possible observability constant 
$C$ in the inequality \eqref{ineqobs}, over $\mathcal{U}_L$, because of the duality between controllability and observability. An interesting problem then consists of maximizing the functional $C_T(\chi_\omega)$ over the set $\mathcal{U}_L$. This problem has been discussed in \cite{PTZobsND,PTZparabND}, and for the same reasons as above it has appeared more relevant to introduce the concept of randomized observability constant $C_{T,\textrm{rand}}(\chi_\omega)$.
At this step, one may think of coming back, by duality, to the controllability problem.
Unfortunately, the problem of maximizing $C_{T,\textrm{rand}}(\chi_\omega)$ does not admit any nice interpretation in terms of controlling, say, almost every initial data to $0$ in time $T$. This is due to the fact that the randomization procedure does not commute with the duality operator realizing the duality between observability and controllability.

More precisely, the Gramian $G_T$ defined above does not commute with the randomization procedure. To describe which kind of initial data can be steered to $0$ in a random way, it would be required to compute the image under $G_T$ of the random laws used in the randomization procedure, and then show that these random laws share appropriate probability properties, as in \cite{BurqTzvetkov1}.

Hence, here, we have found more relevant to combine the randomization procedure with the moment method, in which case the problem of the lack of commutation arising in the HUM procedure disappears.

\paragraph{Use of other biorthogonal families.} We have here used the moment problem approach but in a very special way, taking advantage of the fact that eigenvalues grow sufficiently fast so to ensure the existence of a family of time-biorthogonals that allow to build by separation of variables biorthogonal families for all possible supports of the control $\omega$.\\
Of course the issue can be formulated without that restrictive assumption taking advantage of the existence of biorthogonal families in the $(x, t)$ variables, in other words of families $\Lambda=(\theta_k^{\omega,T})_{k\in\N^*}$ such that
$$
\int_0^T\int_\omega \theta^{\omega,T}_k(t,x)e^{-\lambda_jt}\phi_j(x)\, dxdt=\delta_{jk}.
$$ 
But their dependence with respect to $\omega$ seems to be hard to analyze.

\appendix
\section{Proofs}\label{sec:proofs}
In what follows and similarly to what has been done in Section \ref{sec:mainResIntro}, we will consider a convexified version of the problem \eqref{optDesignContrand} to overcome the difficulty related to the non-compactness of the set $\mathcal{U}_L$ defined by \eqref{defULnew} for the $L^\infty$ weak-star topology. We refer to this section for more comments on this procedure.

The convex closure of $\mathcal{U}_L$ for the weak star topology of $L^\infty$ is
\begin{equation}\label{defULbarnew}
\overline{\mathcal{U}}_L = \left\{ a\in L^\infty(\Omega;[0,1])\ \vert\ \int_{\Omega} a(x)\,dx=L\vert\Omega\vert\right\}.
\end{equation}
Replacing $\chi_\omega\in\mathcal{U}_L$ with $a\in\overline{\mathcal{U}}_L$, we consider the convexified formulation of the problem \eqref{optDesignContrand} given by
\begin{equation}\label{defJanew}
\sup_{a\in \overline{\mathcal{U}}_L} \mathcal{J}(a) ,
\end{equation}
where the functional $J$ is naturally extended to $\overline{\mathcal{U}}_L$ by
\begin{equation}\label{defJrelaxnew}
\mathcal{J}(a)=\inf_{j\in \N^*}\gamma_j(T) \int_{\Omega}a(x)\vert \phi_j(x)\vert^2\, dx,
\end{equation}
for every $a\in \overline{\mathcal{U}}_L$. 
We consider as well a relaxed formulation of the truncated optimal design problem \eqref{pb:truncated} by
\begin{equation}\label{defJaNnew}
\sup_{a\in \overline{\mathcal{U}}_L} \mathcal{J}_N(a) ,
\end{equation}
where the functional $\mathcal{J}_{N}$ is naturally extended to $\overline{\mathcal{U}}_L$ by
\begin{equation}\label{defJaNrelaxnew}
\mathcal{J}_N(a)=\inf_{1\leq j\leq N}\gamma_j(T) \int_{\Omega}a(x)\vert \phi_j(x)\vert^2\, dx,
\end{equation}
for every $a\in \overline{\mathcal{U}}_L$. 

We have the following result.

\begin{lemma}\label{propExistRelax2}
For every $L\in (0,1)$, the relaxed problem \eqref{defJanew} (respectively \eqref{defJaNnew}, for any $N\in\N^*$) has at least one solution $a^*\in\overline{\mathcal{U}}_L$ (respectively $a^{N}\in\overline{\mathcal{U}}_L$).
\end{lemma}

We refer to Section \ref{sec:mainrescomillu} for a proof of this result.

\subsection{Proofs of Lemmas \ref{lemm:contHeat} and \ref{lemma:momentCont}}\label{secproof31}
We prove Lemma \ref{lemma:momentCont}, which is a generalization of Lemma \ref{lemm:contHeat}.

We seek a control function $u=u(t,x)$ in $L^2((0,T)\times \Omega)$ achieving the null controllability for the system \eqref{heatEqcontrolled} with the initial condition $y(0,x)=y_0(x) = \sum_{j=1}^{+\infty} a_j\phi_j(x)$, that is, such that $y(T,\cdot)=0$.
Setting $y(t,x)= \sum_{j=1}^{+\infty} y_j(t)\phi_j(x)$, we get
$$
y_j(T) = e^{-\lambda_jT}a_{j}+\int_0^T\int_{\omega} e^{-\lambda_j(T-t)} u(t,x)\overline{\phi}_j(x) \, dx\, dt,
$$
for every $j\in\N^*$. In order to realize the null controllability in time $T$, the control $u$ must be such that
\begin{equation}\label{eq1}
\int_0^T\int_{\omega} e^{-\lambda_j(T-t)} u(t,x)\overline{\phi}_j(x) \, dx\, dt = -a_j e^{-\lambda_jT},
\end{equation}
for every $j\in\N^*$. In order to solve these equations, assume that there exists a sequence $(\theta_j^T)_{j\in\N^*}$ of functions biorthogonal to the family $\Lambda$, that is,
$$ \int_0^T e^{-\lambda_jt}\theta_k^T(t)\, dt = \delta_{jk},$$
for all $(j,k)\in(\N^*)^2$, where $\delta_{jk}=1$ whenever $j=k$, and $\delta_{jk}=0$ otherwise.
Then the function $u$ defined by \eqref{defu} is a formal solution of the moment problem \eqref{eq1}.

\subsection{Proof of Lemmas \ref{lem2} and \ref{lemmContrand}}\label{sec_proof_lem2}
We randomize the initial datum $y_0(x) = \sum_{j=1}^{+\infty} a_j\phi_j(x)$ according to $y_0^\nu(x) = \sum_{j=1}^{+\infty} \beta_j^\nu a_j\phi_j(x)$. Then, the corresponding control $u=\Gamma_{\omega}(y^\nu_{0})$ coming from the moment method, steering $y_0^\nu$ to $0$ in time $T$, is
$$
u^\nu(t,x) = -\sum_{j=1}^{+\infty} \beta_j^\nu a_j e^{-\lambda_jT} \theta^T_j(T-t)\frac{\phi_j(x)}{\int_\omega\vert\phi_j(y)\vert^2 \, dy},
$$
and hence
\begin{equation*}
\Vert u^\nu\Vert^2_{L^2((0,T)\times\Omega)} 
= \sum_{j,k=1}^{+\infty} \beta_j^\nu \beta_k^\nu a_j \bar{a}_k e^{-(\lambda_{j}+\bar{\lambda}_k)T} \int_0^T \theta_j^T(t)\theta_k^T(t)\, dt \frac{\int_\omega\phi_{ j}(x)\overline{\phi}_{ k}(x)\, dx}{\int_\omega|\phi_j(x)|^2 \, dx \int_\omega|\phi_{k}(x)|^2\, dx},
\end{equation*}
and therefore,
$$
\mathbb{E}\left(\Vert u^\nu\Vert^2_{L^2((0,T)\times\Omega)} \right) = \sum_{j=1}^{+\infty} \frac{\vert a_j\vert^2 e^{-2\Real(\lambda_{j})T} }{\int_\omega|\phi_j(x)|^2 \, dx}\int_0^T \theta_j^T(t)^2\, dt .
$$
The result follows.

\subsection{Proof of Proposition \ref{prop:trunc}}\label{secproof:prop:trunc}
This proof is similar to the one of \cite[Proposition 2]{PTZparabND}. We include it in the present paper, for the sake of completeness and readability.

For every $N\in\N^*$, we consider the relaxed truncated problem \eqref{defJaNrelaxnew}, where the functional $J_N$ is defined by \eqref{defJaNnew}.
Using the same arguments as in the proof of Lemma \ref{propExistRelax2}, it is clear that the problem \eqref{defJaNrelaxnew} has at least one solution $a^N\in\overline{\mathcal{U}}_L$. Let us prove that $a^N$ is the characteristic function of a set $\omega^N$ such that $\chi_{\omega^N}\in\mathcal{U}_L$. Defining the simplex set $\mathcal{S}_N=\left\{\alpha=(\alpha_j)_{1\leq j\leq N}\in [0,1]^N\mid \sum_{j=1}^N\alpha_j=1\right\}$,
it follows from the Sion minimax theorem (see \cite{Sion}) that
\begin{eqnarray*}
\sup_{a\in \overline{\mathcal{U}}_L} \min_{1\leq j\leq N} \gamma_j(T) \int_{\Omega}a(x) \vert\phi_j(x)\vert^2\,dx & = & \max_{a\in \overline{\mathcal{U}}_L} \min_{\alpha\in\mathcal{S}_N}\int_{\Omega}a(x) \sum_{j=1}^N\alpha_j \gamma_j(T) \vert\phi_j(x)\vert^2\,dx\\
 & = &  \min_{\alpha\in\mathcal{S}_N}\max_{a\in \overline{\mathcal{U}}_L}\int_{\Omega}a(x) \sum_{j=1}^N\alpha_j \gamma_j(T) \vert\phi_j(x)\vert^2\,dx ,
\end{eqnarray*}
and that there exists $\alpha^N\in \mathcal{S}_N$ such that $(a^N,\alpha^N)$ is a saddle point of the functional
$$
(a,\alpha)\in \overline{\mathcal{U}}_L\times \mathcal{S}_N\longmapsto\sum_{j=1}^N\alpha_j  \gamma_j(T) \int_{\Omega}a(x)\vert\phi_j(x)\vert^2\,dx.
$$
Therefore, $a^N$ is solution of the optimal design problem
$$
\max_{a\in \overline{\mathcal{U}}_L}\int_{\Omega}a(x) \sum_{j=1}^N\alpha_j^N \gamma_j(T) \vert\phi_j(x)\vert^2\,dx.
$$
We set $\varphi_N(x)=\sum_{j=1}^N\alpha_j^N \gamma_j(T) \vert\phi_j(x)\vert^2$, for every $x\in\Omega$.
It follows from $\Hun$ that $\varphi_N$ is never constant on any subset of $\Omega$ of positive measure. 
Therefore, there exists $\lambda^N$ such that $a^N(x)=1$ whenever $\varphi_N(x)\geq\lambda_N$, and $a^N(x)=0$ otherwise. In other words, $a^N=\chi_{\omega^N}\in\mathcal{U}_L$, with $\omega^N=\{x\in\Omega\ \vert\ \varphi_N(x)>\lambda_N\}$. 

The uniqueness of $a_N$ follows from the fact that, as proved above, any optimal solution is a characteristic function. Indeed if there were two optimal sets, then any convex combination would also be an optimal solution because $J_N$ is concave. This raises a contradiction since any maximizer has to be a characteristic function.

Under the additional assumption $\Htrois$, the function $\varphi_N$ is analytic in $\Omega$ and therefore $\omega^N$ is an open semi-analytic set. 

\subsection{Proof of Theorem \ref{cormaintheo}}\label{sec:proofcormaintheo}
This proof is a immediate adaptation of the one of Theorem \ref{thm:odp1D}.

Indeed, notice that according to the assumption $\Hdeux$, the estimate \eqref{eq:limproof}, constituting the starting point of the proof of Theorem \ref{thm:odp1D}, can be replaced by 
$$
\liminf_{j\to +\infty}\gamma_{j}(T)\int_{\Omega}a^*(x)|\phi_{j}(x)|^2\, dx>\gamma_{1}(T),
$$
where $a^*$ stands for a solution of the relaxed problem \eqref{defJanew}. 

The rest of the proof is then unchanged.

\subsection{Proof of Theorem \ref{thm:lumped}}\label{secproof:lumpedthm}
We expand $g=\sum_{j=1}^{+\infty}g_{j}\phi_{j}$, with $g_{j}=\int_{\Omega} g(x)\overline{\phi}_{j} (x)\, dx$ for every $j\in\N^*$. Note that, since $\Vert g\Vert_{L^2(\Omega)}=1$, we have
$\sum_{j=1}^{+\infty}|g_{j}|^2=1$. 
We define the convex set $\mathcal{S}$ by
$$
\mathcal{S}=\Big\{\beta=(\beta_{j})_{j\in\N\times\N^*}\in \ell^1(\R_{+}) \quad\vert\quad \sum_{j\in \N^*}\beta_{j}=1\Big\} ,
$$
Note that, if $(\eta_{j})_{j\in\N^*}$ is a sequence of positive real numbers, then
$$
\inf_{j\in\N^*}\eta_{j}=\inf_{\beta\in\mathcal{S}}\sum_{j=1}^{+\infty}\beta_{j}\eta_{j}.
$$
Therefore, the optimal value for the problem \eqref{pboptGamma1010} coincides with the optimal value of a convexified problem, as follows. 
Writing $\alpha=(\alpha_{j})_{j\in\N^*}=(g_{j}^2)_{j\in\N^*}$, there holds
\begin{eqnarray*}
\sup_{\Vert g\Vert_{L^2(0,\pi)}=1}\inf_{j\in\N^*}\gamma_{j}(T)\left|\int_{\Omega} g(x)\overline{\phi}_{j} (x)\, dx\right|^2 & = & \sup_{\sum_{j=1}^{+\infty} |g_{j}|^2=1}\inf_{\beta \in\mathcal{S}}\sum_{j=1}^{+\infty}\gamma_{j}(T)\beta_{j}|g_{j}|^2\\
& = & \sup_{\alpha\in\mathcal{S}}\inf_{\beta \in\mathcal{S}}F(\alpha,\beta),
\end{eqnarray*}
where the functional $F$ is defined by
$
F(\alpha,\beta)=\sum_{j=1}^{+\infty}\gamma_{j}(T)\alpha_{j}\beta_{j}.
$
In accordance with \eqref{assump:specGammaj}, define the sequence $\lambda^*=(\lambda_{j}^*)_{j\in\N^*}$ by $\lambda^*_{j}= \left(\sum_{j=1}^{+\infty}\frac{1}{\gamma_{j}(T)}\right)^{-1}\frac{1}{\gamma_{j}(T)}$ for every $j\in\N^*$. Clearly, $\lambda^*\in\mathcal{S}$ and we have
$$
\sup_{\alpha\in\mathcal{S}}\inf_{\beta \in\mathcal{S}}F(\alpha,\beta)\leq \sup_{\alpha\in\mathcal{S}}F(\alpha,\lambda^*)= \left(\sum_{j=1}^{+\infty}\frac{1}{\gamma_{j}(T)}\right)^{-1}.
$$
Similarly, 
$$
\sup_{\alpha\in\mathcal{S}}\inf_{\beta \in\mathcal{S}}F(\alpha,\beta)\geq \inf_{\beta\in\mathcal{S}}F(\lambda^*,\beta)= \left(\sum_{j=1}^{+\infty}\frac{1}{\gamma_{j}(T)}\right)^{-1}.
$$
It follows that
$$
\sup_{\Vert g\Vert_{L^2(0,\pi)}=1}\inf_{j\in\N^*}\gamma_{j}(T)\left|\int_{\Omega} g(x)\phi_{j} (x)\, dx\right|^2= \left(\sum_{j=1}^{+\infty}\frac{1}{\gamma_{j}(T)}\right)^{-1},
$$
and the supremum is reached if, and only if
$$
|g_{j}|^2=\left(\sum_{j=1}^{+\infty}\frac{1}{\gamma_{j}(T)}\right)^{-1}\frac{1}{\gamma_{j}(T)}
$$
for every $j\in\N^*$. The conclusion follows.

\bigskip

\noindent{\bf Acknowledgment.}
The first author was partially supported by the ANR project OPTIFORM.\\
This work was partially supported by  the Advanced Grant DYCON (Dynamic Control) of the European Research Council Executive Agency - GA 694126, ICON of the French ANR-2016-ACHN-0014-01, FA9550-15-1-0027 of AFOSR,  A9550-14-1-0214 of the EOARD-AFOSR, and the MTM2014-52347 Grant of the MINECO (Spain).


\begin{thebibliography}{99}

\small

\bibitem{akhiezer} N.I. Akhiezer,
\textit{Elements of the theory of elliptic functions},
Translations of Mathematical Monographs {\bf 79}, American Mathematical Society (1990).

\bibitem{allaireMunch}
G. Allaire, A. M\"unch, F. Periago,
\textit{Long time behavior of a two-phase optimal design for the heat equation},
SIAM J. Control Optim. {\bf 48} (2010), no. 8, 5333--5356.

\bibitem{ammar}
F. Ammar-Khodja, A. Benabdallah, M. Gonz\'alez-Burgos, L. de Teresa,
\textit{The Kalman condition for the boundary controllability of coupled parabolic systems. Bounds on biorthogonal families to complex matrix exponentials}, 
J. Math. Pures Appl. \textbf{96} (2011), no. 6, 555--590.

\bibitem{AEWZ}
J. Apraiz, L. Escauriaza, G. Wang, C. Zhang,
\textit{Observability inequalities and measurable sets},
to appear in J. Europ. Math. Soc. (2014).

\bibitem{armaoua}
A. Armaoua, M. Demetriou, 
\textit{Optimal actuator/sensor placement for linear parabolic PDEs using spatial $H^2$ norm}, 
Chemical Engineering Science {\bf 61} (2006), 7351--7367.

\bibitem{Avdonin}
S.A. Avdonin, S.A. Ivanov,
\textit{Families of exponentials. The method of moments in controllability problems for distributed parameter systems},
Cambridge University Press, Cambridge, 1995.

\bibitem{BucurButtazzo}
D. Bucur, G. Buttazzo,
\textit{Variational methods in shape optimization problems},
Progress in Nonlinear Differential Equations {\bf 65}, Birkh\"auser Verlag, Basel (2005).

\bibitem{Burq}
N. Burq, 
\textit{Large-time dynamics for the one-dimensional Schr\"odinger equation},
Proc. Roy. Soc. Edinburgh Sect. A. {\bf 141} (2011), no. 2, 227--251.

\bibitem{BurqTzvetkov1}
N. Burq, N. Tzvetkov,
\textit{Random data Cauchy theory for supercritical wave equations. I. Local theory},
Invent. Math. {\bf 173} (2008), no. 3, 449--475.

\bibitem{fattorini1}
H.O. Fattorini, D.L. Russell, 
\textit{Uniform bounds on biorthogonal functions for real exponentials with an application to the control theory of parabolic equations}, 
Quart. J. Appl. Math. {\bf 32} (1974), 45--69.

\bibitem{fattorini2}
H.O. Fattorini, D.L. Russell, 
\textit{Exact controllability theorems for linear parabolic equation in one space dimension}, 
Arch. Ration. Mech. Anal. {\bf 43} (1971), 272--292.

\bibitem{Hansen}
S.W. Hansen,
\textit{Bounds on functions biorthogonal to sets of complex exponentials; control of damped elastic systems},
J. Math. Anal. Appl. {\bf 158} (1991), no. 2, 487--508.

\bibitem{Hardt}
R.M. Hardt,
\textit{Stratification of real analytic mappings and images},
Invent. Math. {\bf 28} (1975).

\bibitem{harris} 
T.J. Harris, J.F. Macgregor, J.D. Wright,
\textit{Optimal sensor location with an application to a packed bed tubular reactor},
AIChE Journal {\bf 26} (1980), no. 6, 910--916.

\bibitem{Hironaka}
H. Hironaka,
\textit{Subanalytic sets}, in: Number Theory, Algebraic Geometry and Commutative Algebra. In honor of Y. Akizuki, Tokyo (1973).

\bibitem{Khapalov}
A.Y. Khapalov,
\textit{Controllability of partial differential equations governed by multiplicative controls}
Lecture Notes in Mathematics, Vol. 1995, 2010, XV, 284 p.

\bibitem{Kumar}
S. Kumar, J.H. Seinfeld,
\textit{Optimal location of measurements for distributed parameter estimation},
IEEE Trans. Autom. Contr. {\bf 23} (1978), 690--698.

\bibitem{Morris}
K. Morris,
\textit{Linear-quadratic optimal actuator location},
IEEE Trans. Automat. Control {\bf 56} (2011), no. 1, 113--124. 

\bibitem{munchZuazua} 
A. M{\"u}nch, E. Zuazua,
\textit{Numerical approximation of null controls for the heat equation : ill-posedness and remedies},
Inv. Problems {\bf 26} (2010), no. 8.

\bibitem{lions2}
J.-L. Lions,
\textit{Exact controllability, stabilizability and perturbations for distributed systems},
SIAM Rev. {\bf 30} (1988), 1--68.

\bibitem{lions}
J.-L. Lions,
\textit{Contr\^olabilit\'e exacte, perturbations et stabilisation de syst\`emes distribu\'es}, Tome 1,
Recherches en Math\'ematiques Appliqu\'ees [Research in Applied Mathematics], Masson (1988).

\bibitem{MicuZuazua}
S. Micu, E. Zuazua,
\textit{Regularity issues for the null-controllability of the linear 1-d heat equation},
Syst. Cont. Letters {\bf 60} (2011), 406--413.

\bibitem{munchPedr} 
A. M{\"u}nch, P. Pedregal,
\textit{Numerical null controllability of semi-linear 1D heat equations: fixed point, least squares and Newton methods},
Math. Control Relat. Fields {\bf 3} (2012) , no. 2, 217--246.

\bibitem{munchHeat} 
A. M{\"u}nch, F. Periago,
\textit{ Optimal distribution of the internal null control for the 1D heat equation},
J. Diff. Equations {\bf 250}, 95--111 (2011).

\bibitem{Pontryagin}
L.S. Pontryagin, V.G. Boltyanskii, R.V. Gamkrelidze, E.F. Mishchenko,
\textit{The mathematical theory of optimal processes},
Interscience Publishers John Wiley \& Sons, Inc. New York-London (1962).

\bibitem{PTZobspb1} 
Y. Privat, E. Tr\'elat, E. Zuazua,
\textit{Complexity and regularity of maximal energy domains for the wave equation with fixed initial data},
Discrete Cont. Dynam. Syst. {\bf 35} (2015), no. 12, 6133--6153.

\bibitem{PTZ_HUM1D} Y. Privat, E. Tr{\'e}lat, E. Zuazua,
\textit{Optimal location of controllers for the one-dimensional wave equation},
Ann.\ Inst.\ H.\ Poincar\'e Anal.\ Non Lin\'eaire {\bf 30} (2013), 1097--1126.

\bibitem{PTZObs1}
Y. Privat, E. Tr\'elat, E. Zuazua,
\textit{Optimal observability of the one-dimensional wave equation},
J. Fourier Anal. Appl. {\bf 19} (2013), no. 3, 514--544.

\bibitem{PTZobsND} 
Y. Privat, E. Tr\'elat, E. Zuazua,
\textit{Optimal observability of the multi-dimensional wave and Schr\"odinger equations in quantum ergodic domains},
J. Eur. Math. Soc. (JEMS) {\bf 18} (2016), no. 5, 1043--1111.

\bibitem{PTZparabND} 
Y. Privat, E. Tr\'elat, E. Zuazua,
\textit{Optimal shape and location of sensors for parabolic equations with random initial data},
Arch. Ration. Mech. Anal. {\bf 216} (2015), no. 3, 921--981.

\bibitem{Russell}
D.L. Russell,
\textit{Controllability and stabilizability theory for linear partial differential equations: recent progress and open questions},
SIAM Rev. {\bf 20} (1978), no. 4, 639--739.

\bibitem{Sigmund}
O. Sigmund, J.S. Jensen,
\textit{Systematic design of phononic band-gap materials and structures by topology optimization},
R. Soc. Lond. Philos. Trans. Ser. A Math. Phys. Eng. Sci. {\bf 361} (2003), no. 1806, 1001--1019. 

\bibitem{Sion}
M. Sion,
\textit{On general minimax theorems},
Pacific J. Math. {\bf 8} (1958), 171--176.

\bibitem{trelat}
E. Tr\'elat,
\textit{Contr\^ole optimal, th\'eorie \& applications (French) [Optimal control, theory and applications]},
Vuibert, Paris, 2005.

\bibitem{TucsnakWeiss}
M. Tucsnak, G. Weiss
\textit{Observation and control for operator semigroups},
Birkh\"auser Advanced Texts: Basler Lehrb\"ucher, Birkh\"auser Verlag, Basel, Switzerland, 2009.

\bibitem{Ucinski}
D. Ucinski, M. Patan,
\textit{Sensor network design for the estimation of spatially distributed processes},
Int. J. Appl. Math. Comput. Sci. {\bf 20} (2010), no. 3, 459--481.

\bibitem{vandeWal}
M. van de Wal, B. Jager,
\textit{A review of methods for input/output selection},
Automatica {\bf 37} (2001), no. 4, 487--510.

\bibitem{vandewouwer}
A. Vande Wouwer, N. Point, S. Porteman, M. Remy,
\textit{An approach to the selection of optimal sensor locations in distributed parameter systems},
Journal of Process Control {\bf 10} (2000), 291--300.

\bibitem{zuazua} E. Zuazua,
\textit{Controllability and Observability of Partial Differential Equations: Some results and open problems},
Handbook of Differential Equations: Evolutionary Equations {\bf 3}, C. M. Dafermos and E. Feireisl eds.,  Elsevier Science (2006), 527--621.

\end{thebibliography}
\end{document}